\RequirePackage{fix-cm}
\documentclass{svjour3}                     
\usepackage[caption=false]{subfig}
\smartqed  
\usepackage{amssymb}
\usepackage[round]{natbib}
\usepackage{graphicx}
\usepackage{float}
\usepackage{mathrsfs}
\usepackage{makecell,multirow,diagbox}
\usepackage{epstopdf}
\usepackage[caption=false]{subfig}
\usepackage{epsfig}
\usepackage{amsmath}

\usepackage{enumitem}

\usepackage{array}

\usepackage{booktabs}
\usepackage{arydshln}

\usepackage{dsfont}

\usepackage{lineno}

\usepackage[ruled]{algorithm2e}

\usepackage[colorlinks = true, linkcolor = blue, urlcolor = blue, citecolor = blue, anchorcolor = blue]{hyperref}

\usepackage{geometry}
\usepackage{color}
\usepackage{xcolor}

\usepackage{lineno}

\begin{document}
\bibliographystyle{spbasic}
\title{Multi-period Asset-liability Management with Reinforcement Learning in a Regime-Switching Market \thanks{This research is supported by National Natural Science Foundation of China (Grant No.11671204), China Scholarship Council (No.202206840089), and Postgraduate Research \& Practice Innovation Program of Jiangsu Province (Project No.KYCX22\_0392).}
}


\author{Zhongqin Gao        \and
        Ping Chen  \and
        Xun Li   \and
        Yan Lv  \and
        Wenhao Zhang
}


\institute{ Zhongqin Gao \at
              Department of Mathematics, Jiangxi University of Finance and Economics, Nanchang, 330032, China\\
              \email{gaozhongqin@jxufe.edu.cn}           
           \and
           Ping Chen \at
            The Faculty of Business and Economics, The University of Melbourne, Parkville 3010, Australia \\
           \email{ping.chen@unimelb.edu.au}
          \and
          Xun Li \at
          Department of Applied Mathematics, The Hong Kong Polytechnic University, Hong Kong, China\\
         \email{li.xun@polyu.edu.hk}
         \and
           Yan Lv \at
              School of Mathematics and Statistics, Nanjing University of Science and Technology, Nanjing 210094, China\\
              \email{lvyan1998@aliyun.com} 
          \and
           Wenhao Zhang  \at
              Department of Applied Mathematics, The Hong Kong Polytechnic University, Hong Kong, China\\
              \email{wen-hao.zhang@connect.polyu.hk}  
}

\date{Received: date / Accepted: date}

\maketitle

\begin{abstract}
This paper explores the mean-variance portfolio selection problem in a multi-period financial market characterized by regime-switching dynamics and uncontrollable liabilities. To address the uncertainty in the decision-making process within the financial market, we incorporate reinforcement learning (RL) techniques. Specifically, the study examines an exploratory mean-variance (EMV) framework where investors aim to minimize risk while maximizing returns under incomplete market information, influenced by shifting economic regimes. The market model includes risk-free and risky assets, with liability dynamics driven by a Markov regime-switching process. To align with real-world scenarios where financial decisions are made over discrete time periods, we adopt a multi-period dynamic model. We present an optimal portfolio strategy derived using RL techniques that adapt to these market conditions. The proposed solution addresses the inherent time inconsistency in classical mean-variance models by integrating a pre-committed strategy formulation. Furthermore, we incorporate partial market observability, employing stochastic filtering techniques to estimate unobservable market states. Numerical simulations and empirical tests on real financial data demonstrate that our method achieves superior returns, lower risk, and faster convergence compared to traditional models. These findings highlight the robustness and adaptability of our RL-based solution in dynamic and complex financial environments.

\end{abstract}

\section{Introduction}
\label{sec1}
Asset-Liability Management (ALM) is a critical financial strategy used by institutions such as banks, insurance companies, and pension funds to manage risks arising from the mismatch between assets (investments) and liabilities (obligations). The primary goal of ALM is to ensure that an institution can meet its present and future financial commitments while maximizing the returns on assets. For example, \cite{hoevenaars2008strategic} explored asset allocation strategies for investors with risky liabilities, addressing inflation and real interest rate risks. Similarly, \cite{cowley2005securitization} discussed how life insurers can apply securitization to life insurance and annuity portfolios to improve liquidity and manage risks. In the context of pension funds, \cite{andonov2017pension} provided a detailed analysis of how factors such as regulatory incentives, fund maturity, and governance influence asset allocation and risk-taking behavior, particularly focusing on U.S. public pension systems.

Liabilities in ALM are often regarded as uncontrollable due to various external factors. These liabilities, such as pension obligations or insurance claims, must be met regardless of market conditions, making them sensitive to external influences such as interest rates and demographic shifts. For example, in public pension funds, the value of liabilities is highly influenced by the discount rate used to calculate their present value. A decrease in interest rates can substantially increase the present value of pension liabilities, creating management challenges (\cite{novy-marx2011public}). Additionally, regulatory frameworks constrain the choice of discount rates, leading to underestimation of liabilities despite actual obligations remaining unchanged (\cite{andonov2017pension}). Factors like changes in life expectancy and policyholder behavior, particularly in life insurance companies, add another layer of unpredictability (\cite{cowley2005securitization}), making liabilities difficult to manage.

Given these uncontrollable liabilities and unpredictable market conditions, portfolio selection under uncertainty has become a cornerstone of modern finance. This paper focuses on the problem of mean-variance portfolio selection in a regime-switching market, incorporating reinforcement learning (RL) to address the uncertainty inherent in the optimal control process.

The mean-variance (MV) framework, pioneered by \cite{markowitz1952portfolio}, remains fundamental to portfolio optimization, aiming to minimize risk (variance) while achieving a targeted return. Despite its computational efficiency, the MV model has limitations, such as reliance on accurate parameter estimation and sensitivity to changing market conditions. Researchers have extended the MV model to account for time-varying conditions, transaction costs, and more complex asset dynamics, including \cite{Li2000}. A persistent challenge, however, is the time inconsistency in multi-period MV problems, where strategies optimal in the initial period may no longer be optimal in subsequent periods (\cite{basak2010dynamic}). Efforts to resolve this issue, such as precommitment strategies and dynamic formulations, still remain sensitive to market parameter estimation.

The literature on mean-variance asset-liability management (MV-ALM) has advanced significantly, addressing the complexities posed by uncontrollable liabilities and stochastic market conditions. \cite{zhang2016cev} studied a MV-ALM problem under the constant elasticity of variance (CEV) process, deriving feedback portfolio strategies through backward stochastic differential equations (BSDEs). \cite{li2018dynamic} extended this work by introducing a dynamic derivative-based investment strategy for MV-ALM with stochastic volatility. \cite{wei2017time} developed a time-consistent strategy for MV-ALM with uncontrollable liabilities using forward-backward stochastic differential equations (FBSDEs), while \cite{zhang2020debt} proposed a debt-ratio-based strategy for incorporating debt management into MV-ALM. \cite{peng2020partial} investigated MV-ALM under partial information, deriving strategies when liability processes are not fully observable. Collectively, these studies address the challenges posed by stochastic liabilities and dynamic market conditions in portfolio optimization.

Regime-switching models, introduced by \cite{hamilton1989new}, provide a way to capture the cyclical behavior of financial markets, where asset returns are influenced by different economic regimes, such as bull and bear markets. These models, which describe non-linear market behavior, have been shown to enhance the robustness of investment strategies, see \cite{ang2002international} and \cite{hardy2001regime}. However, regime-switching adds complexity, as investors must estimate both asset parameters and transition probabilities between regimes, often using filtering techniques to handle partial observability.

The application of regime-switching models in ALM has significantly advanced the understanding of dynamic market conditions. \cite{chen2008continuous} introduced a continuous-time regime-switching model for ALM, demonstrating the impact of regime-switching on the efficient frontier. \cite{chen2011multiperiod} extended the MV-ALM framework to a multi-period setting with regime-switching, providing analytical solutions for portfolio optimization under changing market conditions. \cite{wei2013timeconsistent} addressed time-consistent strategies in a regime-switching environment using the extended Hamilton-Jacobi-Bellman (HJB) equation, while \cite{yao2016dynamic} developed dynamic ALM models with stochastic cash flows, enhancing portfolio management in uncertain markets.

With the rapid advancement of computing power, RL has become a powerful tool for addressing uncertainty in financial markets. Rooted in machine learning, RL focuses on decision-making in environments with incomplete information, enabling the development of optimal strategies through direct interaction with the environment (see \cite{sutton1998reinforcement}). In recent years, RL has been increasingly applied to various financial problems, such as portfolio management. It eliminates the need for explicit model parameter estimation by continuously learning from real-time market data, as demonstrated by \cite{deng2016deep}. A key strength of RL is its ability to balance exploration and exploitation in optimizing investment strategies, dynamically adjusting to new information as it becomes available. Recent advancements, particularly in deep reinforcement learning (DRL), have significantly improved the scalability and performance of these algorithms, especially in high-dimensional and complex environments, as discussed in \cite{liang2018deep}. However, integrating RL with traditional financial models, such as mean-variance optimization, is still an evolving area, particularly when factoring in real-world complexities like regime-switching and uncontrollable liabilities.

RL has proven to be an effective approach in solving dynamic portfolio optimization problems. \cite{Jia2022} proposed an RL algorithm for continuous-time mean-variance portfolio selection, comparing it with DRL methods such as deep deterministic policy gradient (DDPG) and proximal policy optimization (PPO). Their algorithm demonstrates superior performance in achieving mean-variance efficiency in real financial markets. Building on this, \cite{dai2023learning} developed a learning equilibrium mean-variance strategy that addresses time inconsistency in multi-period portfolio selection using RL. Additionally, \cite{Wang2020} framed continuous-time mean-variance portfolio optimization as an RL problem, providing a theoretical foundation and empirical results that highlight the potential of RL in financial decision-making. \cite{hambly2023recent} offered an overview of recent advances in RL applications in finance, underscoring its scalability and adaptability in asset management. Finally, \cite{cui2022multiperiod} presented a comprehensive survey on multi-period mean-variance portfolio selection models, examining the challenges and solutions of applying dynamic programming and RL to multi-period investment problems.

Building on the work of \cite{Cui2023} and \cite{Wu2024}, this paper contributes to the literature by exploring an ALM model with regime-switching and uncontrollable liabilities in a multi-period framework. We develop an RL algorithm to solve the exploratory mean-variance (EMV) portfolio selection problem under both fully and partially observable market regimes. The RL algorithm follows a policy improvement framework, iteratively refining the portfolio selection strategy at each step to ensure convergence toward the optimal solution. To address the constraints of the portfolio optimization problem, we introduce a self-correcting scheme for learning the Lagrange multiplier. The algorithm is designed to operate in both fully observable and partially observable environments. In scenarios with partial information, where market regimes are not fully observable, a stochastic filtering technique is used to estimate the hidden states. This filtering step transforms the problem into a form that enables the RL algorithm to function effectively.

We provide a detailed analysis of the RL-based algorithms through both simulation and empirical studies, and compare three approaches: PoEMV-1 (which incorporates market regime learning), PoEMV-2 (which disregards regime learning), and CoEMV (the classical approach under complete information).

Key findings from the simulation study demonstrate that:

\begin{itemize}
\item The PoEMV-1 algorithm consistently outperforms both PoEMV-2 and CoEMV, achieving better mean returns, lower variance, and superior Sharpe ratios. It exhibits faster convergence and less volatility in terminal wealth.
\item The inclusion of market state learning in PoEMV-1 significantly improves the portfolio's ability to approach the target wealth with less fluctuation compared to PoEMV-2, which does not incorporate this learning process.
\item Empirical studies using real financial data (S\&P 500 index) further validate the effectiveness of PoEMV-1, showcasing its robust performance in regime-switching environments. The PoEMV-1 algorithm adapts to market changes and outperforms traditional control approaches in terms of risk-adjusted returns.
\end{itemize}

The remainder of this paper is organized as follows. Section \ref{sec2} introduces the problem formulation, reviewing the classical multi-period mean-variance portfolio selection model and extending it to the regime-switching exploratory mean-variance (EMV) model with uncontrollable liabilities. In Section \ref{sec3}, we derive the solution to the EMV problem, first under the assumption of complete information and later addressing partial information through stochastic filtering techniques. Section \ref{sec4} presents the reinforcement learning (RL) algorithm designed to solve the EMV problem in both completely and partially observable market regimes. Numerical results and simulations demonstrating the performance of the proposed RL-based portfolio strategies are provided in Section \ref{sec5}. Finally, Section \ref{sec6} concludes the paper, summarizing the contributions and discussing potential future research directions.

\section{Problem formulation}\label{sec2}

Let $T$ denote the number of decision periods, with decision dates $0,\ldots,T-1$ and terminal date $T$. Let $(\Omega, \mathcal{F}, \{\mathcal{F}_t\}_{t=0,\ldots,T}, \mathbb{P})$ be a filtered probability space, where $\mathcal{F}_t$ represents the full market information up to time $t$. The investor's information is specified separately for each observation structure. The investors want to earn more return and take less risk at the planning time $T$. We need to select an objective function that balances the return and the risk, and the mean-variance criterion is an appropriate choice. In this section, we first review the classical multi-period mean-variance (MV) portfolio selection problem, and then introduce the regime-switching exploratory mean-variance (EMV) problem with uncontrollable liability studied in this paper.

\subsection{Classical multi-period MV problem}

In a time horizon with $T$ time periods, we consider a financial market composed of two assets: one risk-free asset and one risky asset. The gross return of the risk-free asset is given by a constant $r_f>0$, and $r_t$ is the excess return of risky asset from period $t$ to $t+1$, which subjects to a distribution with mean $a$ and variance $\sigma^2$, and $\{r_t, t=0, 1, \cdots, T-1\}$ are statistically independent. Denote by $\{x_t^u, t=0,1, \cdots, T\}$ the wealth process of an investor who rebalances the portfolio investing in the risky and riskless assets under investment strategy $u=\{u_0, u_1, \cdots, u_{T-1}\}$, where $u_t$ is the amount invested in the risky asset at time period $t$ and $x_t^u$ is supposed to be adapted to $\mathcal{F}_t$. The allocation $u_t$ is $\mathcal F_t$-measurable, and the excess return $r_t$ over the period from $t$ to $t+1$ is $\mathcal F_{t+1}$-measurable and independent of $\mathcal F_t$.
We assume that transactions are carried out at the beginning of every time period with no transaction cost or tax charged. Then, the dynamics of wealth process is defined as
\begin{align}\label{xt1}
&x_{t+1}^u=r_f x_t^u+r_t u_t,\ \ \ t=0, 1, \cdots, T-1,
\end{align}
with an initial endowment being $x_0^u=x_0>0$.

The classical multi-period MV portfolio selection problem aims to find the optimal admissible investment strategy such that the variance of the terminal wealth is minimized for a given investment target $d$, i.e.
\begin{align}\label{min11}
&\min_{u} \mbox{Var}(x_T^{u})=\mathbb{E}[(x_T^{u})^2]-d^2
\quad \text{s.t.} \quad  \mathbb{E}[x_T^{u}]=d \ \mbox{and} \ \eqref{xt1}.
\end{align}
It is well-known that the problem \eqref{min11} is time-inconsistent due to the variance term in the objective, and the initial precommitted policy need not remain optimal when the criterion is reinitialized at a later date. In this paper, we focus on finding the optimal precommitted strategy, which is optimal only at $t=0$.
Following \cite{Li2000}, we introduce a Lagrange multiplier $w\in\mathbb{R}$ to transform the constrained problem \eqref{min11} into an unconstrained problem as
\begin{align}\label{minuncon}
&\min_{u}\mathbb{E}[{x_T^{u}}^2]-d^2-2w(\mathbb{E}[x_T^{u}]-d)
=\min_{u}\mathbb{E}[(x_T^{u}-w)^2]-(w-d)^2,
\ \ \ \ s.t. \ \eqref{xt1}.
\end{align}
This problem can be solved analytically and its optimal solution denoted by $u^*=\{u_0^*, u_1^*, \cdots, u_{T-1}^*\}$ depends on $w$. The optimal multiplier $w^*$ can be obtained by solving $\mathbb{E}[x_T^{u^*}]=d$.

Given the complete knowledge of the model parameters, the classical MV problem \eqref{min11} and many of its variants have been solved rather completely.
When implementing these solutions, one needs to estimate the market parameters from historical time series of asset prices, a procedure known as identification in classical adaptive control. However, in practice, it is difficult to estimate the investment opportunity parameters, especially the mean return with a workable accuracy.
Moreover, the classical optimal MV strategies are often extremely sensitive to these parameters, largely due to the procedure of inverting ill-conditioned variance-covariance matrices to obtain optimal allocation weights. In view of these two issues, the solutions can be greatly irrelevant to the underlying investment objective. Our focus is on how to solve the unconstrained problem \eqref{minuncon} for a fixed $w$ by reinforcement learning (RL) in a regime-witching market.

\subsection{Regime-switching EMV problem with uncontrollable liability}

RL techniques can use observed or simulated trajectories to learn portfolio allocations without evaluating the full market transition law in each policy update. RL algorithms, driven by historical data, output optimal allocations directly. This is made possible by direct interactions with the unknown investment environment, in a exploring (learning) while exploiting (optimizing) fashion. Under partial information, we use specified or estimated market parameters to construct the regime filter.
To effectively implement the optimal solution, we consider the EMV portfolio selection problem within the framework of RL. Meanwhile, to better simulate the market environment, we incorporate regime-switching and uncontrollable liability into our analysis.
Now, we will proceed to introduce the market state, asset process, liability process and the corresponding EMV problem.

\subsubsection{Financial market}

Suppose two assets are available for investment, labeled as $0$ and $1$, respectively.
From an economic standpoint, the return rates of the assets can be influenced by the market state, while the market state inevitably change over time.

We assume that the gross asset returns depend on the market regime process, which is described by a multi-period homogeneous Markov chain $\varepsilon=\{\varepsilon_{t}\}_{t=0}^T$ with the state space $\mathcal{M}=\{1, 2\}$. The two-regime case is frequently examined in both research and practice, which interprets one regime as a bull market and the other as a bear market.
For the Markov chain $\varepsilon$, we characterize its transition probability matrix as $P=(P_{ij})_{2\times2}, (i,j\in\mathcal{M})$, where $P_{ij}\geq0$
stands for the one-step transition probability of $\varepsilon$ from regime $i$ to regime $j$, and $\sum_{j=1}^2P_{ij}=1$.
The gross returns of assets 0 and 1 over period $t$ are denoted by $e_t(\varepsilon_t)=(e_t^0(\varepsilon_t), e_t^1(\varepsilon_t))^\top$, where $e_t^k(\varepsilon_t)$ is the random return of asset $k\in\{0,1\}$ at time period $t$.
The dependence of $e_t(\varepsilon_t)$ on the market regime can induce serial dependence in asset returns, in contrast to the independent-return setting of \cite{Li2000}.
Conditional on $\varepsilon_t=i$, $\{e_t(i)\}_{t=0}^{T-1}$ are assumed independent, as well as $e_t^0(i)$ and $e_t^1(i)$ being independent. We assume finite second moments for the asset and liability gross returns. Given the full past, the current regime $i$ and the realized action, the joint law of the period-$t$ gross-return vector and $\varepsilon_{t+1}=j$ factors as $K_{t,i}(\mathrm dr)P_{ij}$, where $K_{t,i}$ is the joint gross-return law of the two assets and the liability. The mean and covariance of $e_t(i)$ are respectively denoted as
\[
\mathbb{E}[e_t(i)]=
\begin{pmatrix}
\mathbb{E}[e_t^0(i)] \\
\mathbb{E}[e_t^1(i)] \\
\end{pmatrix},\ \ \ \
\mathrm{cov}(e_t(i))=
\begin{pmatrix}
\sigma_t^{00}(i) & \sigma_t^{01}(i)\\
\sigma_t^{10}(i) & \sigma_t^{11}(i)\\
\end{pmatrix},
\]
where $\sigma_t^{01}(i)=\sigma_t^{10}(i)=0$. Certainly, $\mathbb{E}[e_t(i)e_t'(i)]=\mathrm{cov}(e_t(i))+\mathbb{E}[e_t(i)]\mathbb{E}[e_t'(i)]$. We assume $\mathbb E[e_t(i)e_t(i)^\top]>0$ for every $t$ and $i$.
In the quadratic recursions below, this condition ensures that the coefficient of $u^2$ in the continuation cost and that of $x^2$ in the resulting value function are both strictly positive.

\subsubsection{Surplus process}\label{sec0222}

Under the framework of RL, denote by $u_t^{\pi}$ the amount invested in asset $1$ at the beginning of period $t$, which is randomized and represents exploration (learning), leading to a measure-valued or distributional control process $u^{\pi}=\{u_t^{\pi}\}_{t=0}^{T-1}$ whose density function is expressed as $\pi=\{\pi_t \}_{t=0}^{T-1}$.
Note that the investment strategy $u^{\pi}$ satisfies the self-financing assumption and other standard conditions for admissibility outlined in Definition \ref{def1}.
The amount invested in asset $0$ at the beginning of period $t$ is then equal to $x_t^{\pi}-u_t^{\pi}$.
The theoretical action space is $\mathbb R$, so short positions and borrowing are allowed. Thus, the asset value dynamics can be written as
\begin{align}\label{xt2}
x_{t+1}^{\pi}=e_t^0(\varepsilon_t) (x_t^{\pi}-u_t^{\pi}) +e_t^1(\varepsilon_t)u_t^{\pi}
=e_t^0(\varepsilon_t) x_t^{\pi}+\big(e_t^1(\varepsilon_t)-e_t^0(\varepsilon_t)\big)u_t^{\pi},
\end{align}
where $x_0^{\pi}=x_0$ is the initial asset of investor.
The (distributional) control (or portfolio/strategy) process $\pi$ is to model exploration, whose overall level is in turn captured by its accumulative differential entropy
\begin{align*}
\mathcal H(\pi_t):=-\int_{\mathbb R}\pi_t(u)\ln\pi_t(u)\,\mathrm du,\qquad \mathcal H(\pi):=\sum_{t=0}^{T-1}\mathcal H(\pi_t).
\end{align*}

Incorporating uncontrollable liability into our model allows for a more comprehensive evaluation of financial status, risk factors and decision-making processes. This kind of problem is also referred to as asset-liability management problem.
Let $l_t$ be the liability value of the investor at the beginning of the $t$-th time period, denoted as
\begin{align}\label{lt}
& l_{t+1}=q_t(\varepsilon_t) l_t,
 \end{align}
where $q_t(i)$ is the liability gross-return factor, independent of $e_t(i)$ conditional on $\varepsilon_t=i$. The liability dynamics are exogenous to the portfolio decision. The state space permits signed liability values, as the normal specification used below does not enforce positivity of the gross-return factor. From \eqref{xt2} and \eqref{lt}, the surplus value, representing the difference between asset and liability, can be expressed as
\begin{align}\label{st4}
&S_{t+1}^{\pi}:=x_{t+1}^{\pi}-l_{t+1}
=e_t^0(\varepsilon_t) x_t^{\pi}+\big(e_t^1(\varepsilon_t)-e_t^0(\varepsilon_t)\big)u_t^{\pi}
-q_t(\varepsilon_t) l_t.
\end{align}
The terminal surplus is $S_T^{\pi}=x_{T}^{\pi}-l_T$. The definition of admissible strategies is given below.

Under partial information, the investor observes wealth, liability and past executed allocations. The pre-decision information is
$\mathcal G_t:=\sigma\{x_s^\pi,l_s:0\le s\le t;\ u_s^\pi:0\le s<t\}$.
Under complete information, the investor additionally observes $\varepsilon_0,\ldots,\varepsilon_t$. Write $\mathcal I_t$ for the pre-decision information in the case under consideration. Conditional on $\mathcal I_t$, the action is drawn from $\pi_t$ using randomization independent of the market. The investor observes the action upon execution. The wealth, liability and surplus dynamics apply under either information structure.

\begin{definition}\label{def1}
A strategy $\pi=\{\pi_t\}_{t=0}^{T-1}$ is admissible, if it satisfies the following conditions:
\begin{enumerate}[label=(\roman*)]
\item $\pi_t(\mathrm du)=\pi_t(u)\mathrm du$ is an $\mathcal I_t$-measurable probability kernel on $\mathbb R$;
\item $\mathbb{E}[\sum_{t=0}^{T-1}\int_{\mathbb{R}}u^2\pi_t(u)\mathrm{d}u \mid \mathcal I_t]<\infty$;
\item $\mathbb{E}\big[(S_T^{\pi})^2
+\sum_{t=0}^{T-1}|\int_{\mathbb{R}}\pi_t(u)\ln\pi_t(u)\mathrm{d}u| \mid \mathcal I_t\big]<\infty$, with $0\ln 0=0$.
\end{enumerate}
\end{definition}

\subsubsection{Optimization problem}

In the framework of RL, the precommitted EMV portfolio selection problem with uncontrollable liability and regime-switching market minimizes
\[
\operatorname{Var}(S_T^\pi)+\lambda\,
\mathbb E\big[\sum_{s=0}^{T-1}\int_{\mathbb R}\pi_s(u)\ln\pi_s(u)\,\mathrm du\big]
\quad\text{s.t}\quad \mathbb E[S_T^\pi]=d.
\]
Here $\lambda>0$ is the exploration weight (or temperature parameter) reflecting the trade-off between exploration and exploitation. For any fixed $w\in\mathbb R$, the associated entropy-regularized objective function is
\begin{align}\label{j3}
&J_t^{\pi}(\mathcal I_t; w)
:=\mathbb{E}\big[(S_T^{\pi}-w)^2
+\lambda\sum_{s=t}^{T-1}\int_{\mathbb{R}}\pi_s(u)\ln\pi_s(u)\mathrm{d}u\mid \mathcal I_t \big]-(w-d)^2 \nonumber
\\
&=\mathbb{E}\big[(x_T^{\pi})^2-2wx_T^{\pi}-2x_T^{\pi}l_T+l_T^2+2wl_T
+\lambda\sum_{s=t}^{T-1}\int_{\mathbb{R}}\pi_s(u)\ln\pi_s(u)\mathrm{d}u
\mid \mathcal I_t\big]-d^2+2wd.
\end{align}
At the initial date the expected objective differs from the expected squared deviation from $d$ plus the entropy term by $2(d-w)(\mathbb E[S_T^\pi]-d)$. Consequently, a fixed-$w$ optimizer satisfying the mean constraint also solves the constrained problem. Such a multiplier must be found for the initial law, provided the target can be attained in this way. It is kept fixed in the subsequent dynamic program. For fixed $w$ and a fixed policy before date $t$, all continuation policies are conditioned on the same pre-decision information $\mathcal I_t$.
The EMV problem aims to minimize the objective function \eqref{j3}, this can be formalized as
\begin{align}\label{min8}
&J^*_t(\mathcal I_t; w):=\min_{\pi} J_t^{\pi}(\mathcal I_t; w),
\end{align}
with terminal value $J_T^*(\mathcal I_T;w)=(x-l-w)^2-(w-d)^2$.
Function $J_t^*$ also termed as the value function.
The complete-information state is $(x,l,i)$, and the partially observed problem admits the state $(x,l,\hat p)$ introduced below.

\section{Solution to the EMV problem}\label{sec3}
This section focuses on solving the EMV portfolio selection problem \eqref{min8}.
We first solve it for fixed multiplier under complete information, where the investor observes the market regime. We then use stochastic filtering to reformulate the partially observed problem \eqref{min8} on the belief state and derive its value function and optimal feedback policy by the Bellman principle.

\subsection{Optimal strategy under completely observable environment}\label{subsecopcom}

We assume that the investors can observe the market state $\varepsilon_t$ before acting at period $t$. The policy $\pi^c$ uses the complete-information filtration specified in Section \ref{sec0222}.
In this completely observable situation, the wealth process and the liability process are governed by \eqref{xt2} and \eqref{lt}, respectively, while the complete-information benchmark objective function is defined by
\begin{align*}
&J_t^{\pi^c}(x, l,i; w)
:=\mathbb{E}_{t,x,l,i}\big[(S_T^{\pi^c}-w)^2
+\lambda\sum_{s=t}^{T-1}\int_{\mathbb{R}}\pi^c_s(u)\ln\pi^c_s(u)\mathrm{d}u \big]-(w-d)^2
\\
&=\mathbb{E}_{t,x,l,i}\big[(x_T^{\pi^c})^2-2wx_T^{\pi^c}-2x_T^{\pi^c}l_T+l_T^2+2wl_T
+\lambda\sum_{s=t}^{T-1}\int_{\mathbb{R}}\pi^c_s(u)\ln\pi^c_s(u)\mathrm{d}u
\big]-d^2+2wd,\nonumber
\end{align*}
where $\mathbb{E}_{t,x,l,i}[\cdot]:=\mathbb{E}[\cdot \mid x_t^{\pi^c}=x, l_t=l, \varepsilon_t=i]$, and the entropy regularization term is retained to ensure comparability with the partially observed exploratory formulation. The corresponding EMV problem is
\begin{align}\label{max66}
J_t^{c*}(x, l,i; w):=\min_{\pi^c\in\Pi^c}J_t^{\pi^c}(x, l,i; w),
\end{align}
with
\begin{align}\label{term098}
&J_T^{c*}(x, l, i; w)=(x-l-w)^2-(w-d)^2,
\end{align}
where $\Pi^c$ denotes the policies satisfying Definition \ref{def1} under the complete-information filtration. The conditional value above is written for Markov policies.

Conditional on the current regime $\varepsilon_t = i\in \mathcal{M}$, write $\zeta_{0,t}^{(i)} := e_t^0(i)$, $\zeta_{1, t}^{(i)}:=e_{t}^1(i)-e_{t}^0(i)$, $\zeta_{2, t}^{(i)} = q_{t}(i)$. For $r=0,1,2$ and $(m,n)\in\{(0,1), (0,2), (1,2)\}$, define $A_{r,t}^{(i)}:=\mathbb{E}[\zeta_{r,t}^{(i)}]$, $B_{r,t}^{(i)}:=\mathbb{E}[(\zeta_{r,t}^{(i)})^2]$, $C_{mn,t}^{(i)} = \mathbb{E}[\zeta_{m,t}^{(i)} \zeta_{n,t}^{(i)}]$.
Note that $A_{1,t}^{(i)}=\mathbb{E}[e_t^1(i)]-A_{0,t}^{(i)}$, and since $e_t^0(i)$, $e_t^1(i)$ and $q_t(i)$ are mutually independent under regime $i$, we have $C_{01,t}^{(i)}=A_{0,t}^{(i)} (A_{1,t}^{(i)}+A_{0,t}^{(i)}) - B_{0,t}^{(i)}$, $C_{02,t}^{(i)}=A_{0,t}^{(i)} A_{2,t}^{(i)}$, $C_{12,t}^{(i)}=A_{1,t}^{(i)} A_{2,t}^{(i)}$, and $B_{1,t}^{(i)}=\mathbb{E}[(e_t^1(i))^2]-2A_{0,t}^{(i)}(A_{0,t}^{(i)}+A_{1,t}^{(i)})+B_{0,t}^{(i)}$. Without loss of generality, assume that for any time-dependent function $y$, $\prod_{k=t}^{T} y_k=1$ and $\sum_{k=t}^{T} y_k=0$ if $t>T$. The following theorem provides the value function and optimal investment strategy for the completely observable EMV (CoEMV) case. Its proof is deferred to \ref{appen11}.

\begin{theorem}\label{thm32a}
For $(t,x,l,i)\in\{0,\ldots,T\}\times\mathbb R^2\times\mathcal M$, the value function $J_{t}^{c*}(x,l,i; w)$ of the CoEMV portfolio selection problem \eqref{max66} with fixed $w\in\mathbb{R}$ is given by
\begin{align}\label{eq:general-value-final0}
J_t^{c*}(x,l,i;w)
&=
\frac{\lambda}{2}\bar L_t^{(i)}
+\bar F_{1,t}^{(i)}x^2
-2\bigl(w\bar F_{2,t}^{(i)}
+l\bar F_{3,t}^{(i)}\bigr)x
-\bar D_{1,t}^{(i)}w^2
+2\bigl(\bar A_{2,t}^{(i)}
-\bar D_{2,t}^{(i)}\bigr)wl
+\bigl(\bar B_{2,t}^{(i)}
-\bar D_{3,t}^{(i)}\bigr)l^2 \nonumber \\
& \quad -d^2+2wd,
\end{align}
where,
\begin{align}
&F_{1, t}^{(i)}:={B}_{0,t}^{(i)}{B}_{1,t}^{(i)} - ({C}_{01,t}^{(i)})^2,
\quad
{F}_{2,t}^{(i)}:={A}_{0,t}^{(i)}  {B}_{1,t}^{(i)} - {A}_{1,t}^{(i)} {C}_{01,t}^{(i)},
\quad
{F}_{3,t}^{(i)} :={B}_{1,t}^{(i)}{C}_{02,t}^{(i)} - {C}_{01,t}^{(i)} {C}_{12,t}^{(i)},
\label{f123t} \\
& \bar F_{r,t}^{(i)}
:=
\frac{F_{r,t}^{(i)}}{B_{1,t}^{(i)}}
\sum_{j\in\mathcal M}P_{ij}\bar F_{r,t+1}^{(j)}, (r=1,2, 3),\quad
\bar L_t^{(i)}
:=\ln
\frac{
B_{1,t}^{(i)}
\sum_{j\in\mathcal M}P_{ij}\bar F_{1,t+1}^{(j)}
}{\pi\lambda}
+
\sum_{j\in\mathcal M}P_{ij}\bar L_{t+1}^{(j)}, \label{eq:barF-general} \\
& \bar A_{2,t}^{(i)}
:=
A_{2,t}^{(i)}
\sum_{j\in\mathcal M}P_{ij}\bar A_{2,t+1}^{(j)}, \quad
\bar B_{2,t}^{(i)}
:= B_{2,t}^{(i)}
\sum_{j\in\mathcal M}P_{ij}\bar B_{2,t+1}^{(j)}, \label{eq:barA-general}\\
&\bar D_{1,t}^{(i)}
:=
\sum_{j\in\mathcal M}P_{ij}\bar D_{1,t+1}^{(j)}
+
\frac{
\left(
A_{1,t}^{(i)}
\sum_{j\in\mathcal M}P_{ij}\bar F_{2,t+1}^{(j)}
\right)^2
}{
B_{1,t}^{(i)}
\sum_{j\in\mathcal M}P_{ij}\bar F_{1,t+1}^{(j)}
}, \label{eq:barD1-general} \\
& \bar D_{2,t}^{(i)}
:=
A_{2,t}^{(i)}
\sum_{j\in\mathcal M}P_{ij}\bar D_{2,t+1}^{(j)}
+
\frac{
A_{1,t}^{(i)}C_{12,t}^{(i)}
\left(
\sum_{j\in\mathcal M}P_{ij}\bar F_{2,t+1}^{(j)}
\right)
\left(
\sum_{j\in\mathcal M}P_{ij}\bar F_{3,t+1}^{(j)}
\right)
}{
B_{1,t}^{(i)}
\sum_{j\in\mathcal M}P_{ij}\bar F_{1,t+1}^{(j)}}, \label{eq:barD2-general}\\
&\bar D_{3,t}^{(i)}
:=
B_{2,t}^{(i)}
\sum_{j\in\mathcal M}P_{ij}\bar D_{3,t+1}^{(j)}
+
\frac{
\left(
C_{12,t}^{(i)}
\sum_{j\in\mathcal M}P_{ij}\bar F_{3,t+1}^{(j)}
\right)^2
}{
B_{1,t}^{(i)}
\sum_{j\in\mathcal M}P_{ij}\bar F_{1,t+1}^{(j)}
}, \label{eq:barD3-general}
\end{align}
with the terminal quantities
\begin{align*}
& \bar F_{1,T}^{(i)}=\bar F_{2,T}^{(i)}=\bar F_{3,T}^{(i)}=1, \ \
\bar L_T^{(i)}=0, \ \
\bar A_{2,T}^{(i)}=1, \ \
\bar B_{2,T}^{(i)}=1, \ \
\bar D_{1,T}^{(i)}=\bar D_{2,T}^{(i)}=\bar D_{3,T}^{(i)}=0.
\end{align*}
The optimal feedback control is Gaussian with density function
\begin{equation} \label{pi32a}
\pi_t^{c*}(u;x,l,i,w)
=
\mathcal N \Big(
u \Big|
-\frac{\mu_t^{(i)}(x,l;w)}{\bar B_{1,t}^{(i)}},
\frac{\lambda}{2\bar B_{1,t}^{(i)}}
\Big),
\end{equation}
where
\begin{align}
&\bar B_{1,t}^{(i)}
:= B_{1,t}^{(i)}
\sum_{j\in\mathcal M}P_{ij}\bar F_{1,t+1}^{(j)}, \label{eq:general-b01} \\
& \mu_t^{(i)}(x,l;w)
:=
C_{01,t}^{(i)}x
\sum_{j\in\mathcal M}P_{ij}\bar F_{1,t+1}^{(j)}
-
wA_{1,t}^{(i)}
\sum_{j\in\mathcal M}P_{ij}\bar F_{2,t+1}^{(j)}
-
lC_{12,t}^{(i)}
\sum_{j\in\mathcal M}P_{ij}\bar F_{3,t+1}^{(j)}. \label{eq:general-mu}
\end{align}

\end{theorem}

\begin{remark}\label{rem209}
Suppose that changes in the market state do not affect the return rates of the assets and liabilities.
Then, by redefining $e_t(\varepsilon_t)$ and $q_t(\varepsilon_t)$ as independent random variables $e_t$ and $q_t$, and both of which are independent of $\varepsilon_t$,
the optimal strategy $\pi_t^{c*}$ given by \eqref{pi32a} can be simplified to that of the EMV problem considering liability but without regime-switching.

\end{remark}

\begin{remark}\label{rem2009}
Let $e_t^0(\varepsilon_t)\equiv r_f$ represent the constant return rate for the risk-free asset, $e_t^1(\varepsilon_t)-e_t^0(\varepsilon_t)\equiv r_t$ denote the excess return of risky asset following a normal distribution with mean $a$ and variance $\sigma^2$, and $l_0=0$ indicate the absence of liabilities. Then, the optimal strategy $\pi_t^{c*}$ given by \eqref{pi32a} can be reduced to that of the EMV problem without liability and regime-switching, i.e.
\begin{align*}
\pi_{t}^{c*}(u; x, w)
&=\mathcal{N}\Big(u\Big{|}
 -\frac{a r_f(x-\rho_t w)} {a^2+\sigma^2},
\frac{\lambda}{2 (a^2+\sigma^2)}\big(\frac{a^2+\sigma^2}{\sigma^2 r_f^2}\big)^{T-t-1}\Big),
\end{align*}
where $\rho_t=(r_f^{-1})^{T-t}$, consistent with Theorem 1 in \cite{Cui2023}.

\end{remark}

\subsection{Stochastic filtering}

To solve the partially observable EMV (PoEMV) portfolio selection problem \eqref{min8}, we use the filtering technique to
filter the unobservable market states from the observable information.
For $t\in\{0,\cdots,T\}$, the filtering process for the Markov chain $\varepsilon$ is defined as
\begin{align}\label{p1t}
\hat{p}_t:=\mathbb{P}(\varepsilon_t=1 \mid \mathcal G_t),
\end{align}
which is the probability of the market being in regime $1$ given the observable partial information up to time $t$. Here, $\hat{p}_0:=\mathbb{P}(\varepsilon_0=1)$ is the initial distribution of the market state and is assumed to be known.

Assume $\hat p_0\in(0,1)$ and $P_{11},P_{21}\in(0,1)$ so that the predicted belief $\hat{p}_t$ lies in $(0,1)$. This avoids degenerate beliefs at the endpoints and is consistent with the discrete-time filtering recursion, provided the likelihood terms remain strictly positive. Since the action is generated using independent randomization, $\mathbb P(\varepsilon_t=1\mid\mathcal G_t,u_t^\pi)=\hat p_t$. The one-step filter update admits the standard \emph{Bayes update-prediction} decomposition (the discrete analogue of the Wonham filter for finite-state hidden Markov models).
More precisely, since the new observation $(x^\pi_{t+1},l_{t+1})$ is generated conditional on the current regime $\varepsilon_t$ through \eqref{xt2}-\eqref{lt}, the update from $t$ to $t+1$ can be written as two consecutive steps:
\begin{itemize}
\item \emph{(Bayes update).} The next observation $Y_{t+1}=(x_{t+1}^\pi,l_{t+1})$ updates the belief about $\varepsilon_t$ conditional on the executed action (from $\mathbb P(\varepsilon_t\mid \mathcal G_t)$ to $\mathbb P(\varepsilon_t\mid \mathcal G_{t+1})$).
\item \emph{(Prediction).} The updated probabilities are propagated through the transition matrix $P=(P_{ij})$ to obtain the belief about $\varepsilon_{t+1}$, $\mathbb P(\varepsilon_{t+1}\mid \mathcal G_{t+1})$.
\end{itemize}
This yields a finite-dimensional recursive filter for $\hat p_t$; see, e.g., \cite{Elliott1997}. In particular, we can deduce the following lemma.

\begin{lemma}[Discrete-time filtering recursion]\label{p2t}
Let $G_t^{(i)}(\mathrm dy;x,l,u)$ denote the law of $Y_{t+1}$ given the pre-decision history, $\varepsilon_t=i$, and $u_t^\pi=u$, evaluated at $(x_t^\pi,l_t)=(x,l)$. Choose jointly measurable Radon--Nikodym densities $g_t^{(i)}>0$ with respect to the dominating kernel $G_t^{(1)}+G_t^{(2)}$.
Then the filtering process $\hat p_t$ admits the following two-step recursion:
\begin{align}
\hat p_t^{\,\mathrm{post}}
&:=\mathbb P(\varepsilon_t=1\mid \mathcal G_{t+1})
= \frac{
g_{t}^{(1)}(Y_{t+1};x_t^\pi,l_t, u_t^{\pi})\,\hat p_t
}{
g_{t}^{(1)}(Y_{t+1};x_t^\pi,l_t,u_t^{\pi})\,\hat p_t
+
g_{t}^{(2)}(Y_{t+1};x_t^\pi,l_t, u_t^{\pi})\,(1-\hat p_t)
}, \label{eq:bayes-update} \\
\hat p_{t+1}
&=\mathbb P(\varepsilon_{t+1}=1\mid \mathcal G_{t+1})
= P_{21} + (P_{11}-P_{21})\,\hat p_t^{\,\mathrm{post}}.\label{titer}
\end{align}
In particular, $\hat p_t$ is $\mathcal G_t$-measurable for each $t$.
\end{lemma}

\begin{proof}
Fix $t\in\{0,\ldots,T-1\}$. Conditional on $(\mathcal G_t,u_t^\pi)$, the prior regime probabilities remain $(\hat p_t,1-\hat p_t)$ and the observation law in regime $i$ is $G_t^{(i)}$. Bayes' formula therefore gives \eqref{eq:bayes-update}. In particular, the formula covers observations with zero likelihood in one regime. The joint-kernel assumption gives $\mathbb P(\varepsilon_{t+1}=1\mid\varepsilon_t=i,\mathcal G_{t+1})=P_{i1}$. The law of total probability then gives
\[
\hat p_{t+1}=\sum_{i\in\mathcal M}P_{i1}\mathbb P(\varepsilon_t=i\mid\mathcal G_{t+1})
=P_{21}+(P_{11}-P_{21})\,\hat p_t^{\,\mathrm{post}}.
\]
Since the observation and the executed action are both in $\mathcal G_{t+1}$, this yields the stated measurability.
\end{proof}

By the separation principle, the PoEMV problem \eqref{min8} can be decomposed into a filtering step and a complete observed control problem with belief state $(x_t, l_t, \hat{p}_t)$. The transition of this state can depend on the executed action through the observation law and the Bayes update.
The validity of the separation principle for the MV problem with partial information is established in \cite{Xiong2007} for general price dynamics with unobservable drift. For a Markov feedback policy, the objective \eqref{j3} can be written as
\begin{align*}
&\hat{J}_t^{\pi}(x,l, \hat{p}; w)
:=\mathbb{E}_{t,x,l,\hat{p}}\big[({S}_T^{\pi}-w)^2
+\lambda\sum_{s=t}^{T-1}\int_{\mathbb{R}}\pi_s(u)\ln\pi_s(u)\mathrm{d}u \big]-(w-d)^2\\
&=\mathbb{E}_{t,x,l,\hat{p}}\big[({x}_T^{\pi})^2-2w{x}_T^{\pi}-2{x}_T^{\pi}{l}_T
+{l}_T^2+2w{l}_T
+\lambda\sum_{s=t}^{T-1}\int_{\mathbb{R}}\pi_s(u)\ln\pi_s(u)\mathrm{d}u
\big]-d^2+2wd,\nonumber
\end{align*}
where $\mathbb{E}_{t,x,l,\hat{p}}[\cdot]:=\mathbb{E}[\cdot \mid x_t^{\pi}=x, l_t=l, \hat{p}_t=\hat{p}]$.
The PoEMV problem \eqref{min8} can be reformulated as
\begin{align}\label{min22}
&\hat{J}_t^{*}(x,l,\hat{p}; w):=\min_{\pi\in\Pi}\hat{J}_t^{\pi}(x,l, \hat{p}; w),
\end{align}
with
\begin{align}\label{termcon}
&\hat{J}_T^*(x, l, \hat{p}; w)=(x-l-w)^2-(w-d)^2,
\end{align}
where $\Pi$ denotes the admissible feedback policies under the partial-information filtration in Definition \ref{def1}.

\subsection{Optimal strategy under partially observable environment}\label{subsecparobenv}

Now, the Bellman's principle is applied to solve problem \eqref{min22}, yielding the following result, the proof of which is provided in \ref{appen1010}.

\begin{theorem}\label{them1122}
For $(t,x,l,\hat p)\in\{0,\ldots,T\}\times\mathbb R^2\times(0,1)$, the value function $\hat{J}_{t}^{*}(x,l,\hat{p}; w)$ of the PoEMV portfolio selection problem \eqref{min22} with fixed $w\in\mathbb{R}$ is given by, for $t=T$, $\hat J_T^*(x,l,\hat p;w) = (x-l-w)^2-(w-d)^2$, and for $t=T-1,\ldots,0$, the Bellman value is
\begin{align}\label{hatjt11}
&\hat J_t^*(x,l,\hat p;w)
=
-\lambda\ln
\int_{\mathbb{R}}
\exp\left(
-\mathcal Q_t(u;x,l,\hat p, w)/\lambda
\right)\,\mathrm{d} u,
\end{align}
where
\begin{align}\label{qq22}
\mathcal Q_t(u;x,l,\hat p, w)
&:=
\mathbb{E}\left[
\hat J_{t+1}^*(x_{t+1}^\pi,l_{t+1},\hat p_{t+1};w)
\;\middle|\;
x_t^\pi=x,l_t=l, \hat p_t=\hat p, u_t^\pi=u
\right] \nonumber\\
& = \int_{\mathbb{R}^2}
\hat J_{t+1}^*\bigl(y_1,y_2,\Psi_t(\hat p,y;x,l, u);w\bigr)\,
\widehat G_t(\mathrm{d}y;x,l,\hat p, u),
\end{align}
with $y=(y_1,y_2)\in\mathbb{R}^2$, and the Bayes update and prediction map $\Psi_t$ is given by
\begin{align}\label{psihatpxlu}
&\Psi_t(\hat p,y;x,l, u)
:=
P_{21}
+
(P_{11}-P_{21})
\frac{
g_t^{(1)}(y;x,l, u)\hat p
}{
g_t^{(1)}(y;x,l, u)\hat p
+
g_t^{(2)}(y;x,l, u)(1-\hat p)
}.
\mathfrak{}\end{align}
The map may be set to $\hat p$ when its denominator vanishes, since those observations have zero predictive probability. The posterior-averaged observation kernel $\widehat G_t$ is given by
\begin{align}\label{hatgxlp}
&\widehat G_t(\mathrm{d}y;x,l,\hat p, u)
:=
\hat p\,G_t^{(1)}(\mathrm{d}y;x,l, u)
+
(1-\hat p)\,G_t^{(2)}(\mathrm{d}y;x,l,u).
\end{align}
The normalizing constant is $Z_t(x,l,\hat p;w)
:=
\int_{\mathbb{R}}
\exp\, (
-\mathcal Q_t(v;x,l,\hat p, w)/\lambda)\,\mathrm{d} v$. Under the stated moment and positive-definiteness assumptions, $0<Z_t(x,l,\hat p;w)<\infty$ and the Gibbs feedback is admissible. Then the unique optimal feedback control is Gibbs form with density function
\begin{align}\label{pittt}
\pi_t^*(u; x,l,\hat p, w)
=
\frac{
\exp\big(
-\mathcal Q_t(u;x,l,\hat p, w)/\lambda
\big)
}{
\displaystyle
\int_{\mathbb{R}}
\exp\big(
-\mathcal Q_t(v;x,l,\hat p, w)/\lambda
\big)\,\mathrm{d}v
}.
\end{align}

\end{theorem}
It follows from \eqref{pittt} that, the optimal strategy under partial information depends on the posterior belief $\hat{p}$ about the market state.

\begin{corollary}\label{coroparavalue}
If, at a fixed $(t,x,l,\hat p,w)$, the continuation cost is a positive-curvature quadratic in $u$, write
\[
\mathcal Q_t(u;x,l,\hat p, w)
=
\alpha_t(x,l,\hat p;w)u^2
+
2\beta_t(x,l,\hat p;w)u
+
\gamma_t(x,l,\hat p;w),
\]
with $\alpha_t(x,l,\hat p;w)>0$. Then, the optimal feedback control \eqref{pittt} reduces to the Gaussian density
\begin{align}\label{oppoalp}
\pi_t^*(u; x,l,\hat p, w)
=
\mathcal N\left(
u\;\middle|\;
-\frac{\beta_t(x,l,\hat p;w)}{\alpha_t(x,l,\hat p;w)},
\frac{\lambda}{2\alpha_t(x,l,\hat p;w)}
\right),
\end{align}
and the corresponding value function is
\[
\hat J_t^*(x,l,\hat p;w)
=
\gamma_t(x,l,\hat p;w)
-
\frac{\beta_t(x,l,\hat p;w)^2}{\alpha_t(x,l,\hat p;w)}
+
\frac{\lambda}{2}\ln\frac{\alpha_t(x,l,\hat p;w)}{\pi\lambda}.
\]
Moreover, the variance of the optimal Gaussian policy, reflecting the level of exploration, is decreasing in $t$ whenever $\alpha_t(x,l,\hat p;w)$ is increasing in $t$; the expectation of the optimal Gaussian policy, which reflect the level of exploitation, is independent of the exploration weight $\lambda$ whenever $\alpha_t(x,l,\hat p;w)$ and $\beta_t(x,l,\hat p;w)$ are independent of $\lambda$.
\end{corollary}
This quadratic form holds at the last decision date. At earlier dates it requires a closure condition, since the belief $\Psi_t(\hat p,Y_{t+1};x,l,u)$ can depend on the action.

In addition, to construct a no-learning benchmark $\pi^{sub*}$, we use the unconditional regime probabilities (observation-free prior belief) rather than updating the regime belief $\hat p_t$ from the investor's observations,
\[
\tilde p_t:=\mathbb P(\varepsilon_t=1), \qquad t=0,1,\cdots,T.
\]
Then $\tilde p_t\in(0,1)$ and evolves deterministically according to
\begin{align}\label{tildep00}
&\tilde p_{t+1}=P_{21}+(P_{11}-P_{21})\tilde p_t
= P_{21} \sum_{k=0}^t (P_{11}- P_{21})^k + \tilde{p}_0 (P_{11}- P_{21})^{t+1},
\qquad \tilde p_0=\hat p_0.
\end{align}

In the auxiliary market used to construct this benchmark, returns are drawn independently over time from the date-specific prior mixture of the joint gross-return laws,
\begin{align}\label{tildeassetre}
\widetilde K_t(\mathrm dr)
:=\tilde p_t K_{t,1}(\mathrm dr)+(1-\tilde p_t)K_{t,2}(\mathrm dr).
\end{align}
The feedback obtained in this auxiliary market is evaluated under the original regime-switching dynamics in the simulation study. Under those dynamics, the observed wealth and liability can themselves convey information about the regime.

For each $i\in \mathcal{M}$, $r\in\{0,1,2\}$ and $(m, n)\in \{(0,1), (0,2), (1,2)\}$, let $A_{r,t}^{(i)}$, $B_{r,t}^{(i)}$ and $C_{mn, t}^{(i)}$ be defined exactly as in Section \ref{subsecopcom} under the regime $\varepsilon_t=i$. The moments under \eqref{tildeassetre} are the prior-weighted averages: $\tilde{A}_{r,t}^{(\tilde p)} := \tilde p_t A_{r, t}^{(1)} + (1-\tilde p_t) A_{r, t}^{(2)}$, $\tilde{B}_{r,t}^{(\tilde p)} := \tilde p_t B_{r, t}^{(1)} + (1-\tilde p_t) B_{r, t}^{(2)}$ and $\tilde{C}_{mn,t}^{(\tilde p)} = \tilde p_t C_{mn, t}^{(1)} + (1-\tilde p_t) C_{mn, t}^{(2)}$.
A similar reasoning as in Theorem \ref{thm32a} leads to the following result, the proof is provided in \ref{appen1122}.

\begin{theorem}\label{negtilpthm}
For fixed $w$, the optimal feedback in the auxiliary market \eqref{tildeassetre} is Gaussian with density
\begin{align}\label{pi316}
\pi_t^{sub*}(u; x, l, \tilde p, w)
=
\mathcal N \Big(
u \Big|
-\frac{\tilde{\mu}_t^{(\tilde p)}(x,l;w)}{\tilde{\bar B}_{1,t}^{(\tilde p)}},
\frac{\lambda}{2\tilde{\bar B}_{1,t}^{(\tilde p)}}
\Big),
\end{align}
where
\begin{align}
&\tilde{F}_{1, t}^{(\tilde p)}:={\tilde{B}}_{0,t}^{(\tilde p)}{\tilde{B}}_{1,t}^{(\tilde p)} - ({\tilde{C}}_{01,t}^{(\tilde p)})^2, \quad
{\tilde{F}}_{2,t}^{(\tilde p)}:={\tilde{A}}_{0,t}^{(\tilde p)}  {\tilde{B}}_{1,t}^{(\tilde p)} - {\tilde{A}}_{1,t}^{(\tilde p)} {\tilde{C}}_{01,t}^{(\tilde p)}, \quad
{\tilde{F}}_{3,t}^{(\tilde p)} :={\tilde{B}}_{1,t}^{(\tilde p)}{\tilde{C}}_{02,t}^{(\tilde p)} - {\tilde{C}}_{01,t}^{(\tilde p)} {\tilde{C}}_{12,t}^{(\tilde p)}, \label{tildef123} \\
& \tilde{\bar F}_{r,t}^{(\tilde p)}
:=
\frac{\tilde F_{r,t}^{(\tilde p)}}{\tilde B_{1,t}^{(\tilde p)}} \tilde{\bar F}_{r,t+1}^{(\tilde p)}, (r=1,2, 3),\quad
\tilde{\bar B}_{1,t}^{(\tilde p)}
:= \tilde{B}_{1,t}^{(\tilde p)}
\tilde{\bar F}_{1,t+1}^{(\tilde p)}, \label{eq:general-b01} \\
& \tilde{\mu}_t^{(\tilde p)}(x,l;w)
:=
\tilde C_{01,t}^{(\tilde p)}x
\tilde{\bar F}_{1,t+1}^{(\tilde p)}
-
w \tilde{A}_{1,t}^{(\tilde p)} \tilde{\bar F}_{2,t+1}^{(\tilde p)}
-
l \tilde C_{12,t}^{(\tilde p)} \tilde{\bar F}_{3,t+1}^{(\tilde p)}, \label{eq:general-mu}
\end{align}
with the terminal quantities $\tilde{\bar F}_{1,T}^{(\tilde p)}=\tilde{\bar F}_{2,T}^{(\tilde p)}
=\tilde{\bar F}_{3,T}^{(\tilde p)}=1$.

\end{theorem}

The complete-information policy \eqref{pi32a} uses the observed regime, the partially observed Gibbs policy \eqref{pittt} uses the updated belief, and the no-learning reference \eqref{pi316} uses the deterministic prior. Section \ref{sec5} compares numerical policies associated with these three information structures.

\section{RL algorithm design}\label{sec4}
In this section, we design an RL algorithm to learn the solution and output portfolio allocation strategy in both completely observable and partially observable cases. Beginning with the establishment of the policy improvement theorem alongside its convergence implication for exact Gibbs updates, we subsequently devise an actor-critic algorithm designed to learn optimal solution, incorporating a self-correcting scheme for learning the correct Lagrange multiplier $w$.

\subsection{Policy improvement and policy convergence}\label{subsec 4.1}
A policy improvement scheme can be implemented to update the current policy in a direction that improves the objective function. In other words, for a minimization problem, the policy improvement theorem ensures that the iterated objective function is non-increasing and ultimately converges to the value function. We now present the policy improvement theorems for the multi-period EMV portfolio selection problem.

To solve the complete-information problem, the partially observed problem, and the no-learning benchmark in a unified framework, let $m\in\{c, p, s\}$, where $c$ and $p$ denote complete and partial information and $s$ denotes the auxiliary prior-mixture market. Denote $z_t^c:=\varepsilon_t\in\mathcal M, z_t^p:=\hat p_t\in(0, 1), z_t^s:=\tilde p_t\in(0, 1)$, and write the corresponding state variable as $\xi_t^m:=(x_t^\pi,l_t,z_t^m)$. For the case $m$, write $J_t^{m, \pi}(\xi;w)$ to denote the objective function associated with an admissible policy $\pi$ and $\xi=\xi_t^m$, and let $J_t^{m,*}(\xi;w):=\min_{\pi\in\Pi^m}J_t^{m, \pi}(\xi;w)$, where $\Pi^m$ is the admissible Markov feedback class for model $m$, with the common terminal value $J_T^{m, \pi}(\xi;w):=(x_T^{\pi}-l_T-w)^2-(w-d)^2$.
For any admissible Markov policy $\pi=\{\pi_t\}_{t=0}^{T-1}$, define
\[
J_t^{m, \pi}(\xi;w)
=
(\mathcal T_t^{m, \pi}J_{t+1}^{m, \pi})(\xi)
:=
\int_{\mathbb R}
\big(
Q_t^{m, \pi}(u;\xi,w)+\lambda\ln\pi_t(u\mid \xi)
\big) \pi_t(u\mid \xi)\, \mathrm du,
\]
where $\mathcal T_t^{m, \pi}$ is the one-step evaluation operator, $Q_t^{m, \pi}(u;\xi,w) := \mathbb E [J_{t+1}^{m, \pi}(\xi_{t+1}^m; w) \;|\;\xi_t^m=\xi,\ u_t^{\pi}=u]$.
In particular, the associated state-action values are given by,
\begin{enumerate}[label=(\roman*)]
\item in complete information, $\xi=(x,l,i)$ and $Q_t^{c, \pi}(u;x,l,i,w)
=
\sum_{j\in\mathcal M}
\int_{\mathbb R^2}
J_{t+1}^{c,\pi}(y_1,y_2,j;w)\,
P_{ij}\,G_t^{(i)}(\mathrm dy;x,l, u)$;
\item in partial information, $\xi=(x,l,\hat p)$ and
$Q_t^{p, \pi}(u;x,l,\hat p,w)
=
\int_{\mathbb R^2}
J_{t+1}^{p,\pi}\bigl(y_1,y_2,\Psi_t(\hat p,y;x,l,u); w\bigr)\,
\widehat G_t(\mathrm dy;x,l,\hat p,u)$, where $\Psi_t$ and $\widehat G_t$ are given by \eqref{psihatpxlu} and \eqref{hatgxlp}, respectively;
\item in the auxiliary no-learning benchmark, $\xi = (x,l, \tilde p)$ and
$Q_t^{s, \pi}(u;x,l,\tilde p, w)
=
\int_{\mathbb R^2}
J_{t+1}^{s, \pi}(y_1,y_2, \widetilde\Psi_t(\tilde p); w)\,
\widetilde G_t(\mathrm dy;x,l,\tilde p, u)$,
where $\widetilde\Psi_t(\tilde p):=P_{21}+(P_{11}-P_{21})\tilde p$ and $\widetilde G_t(\mathrm dy;x,l, \tilde p, u)
:=
\tilde p \,G_t^{(1)}(\mathrm dy;x,l,u)
+
(1-\tilde p)\,G_t^{(2)}(\mathrm dy;x,l,u)$.
\end{enumerate}

For the case $m$, the policy-improvement theorem is presented below, with the proof provided in \ref{thm:policy_improvement_unified00}.
\begin{theorem}[Policy improvement]\label{thm:policy_improvement_unified}
For each $t=0,1,\cdots,T-1$ and $\xi=\xi_t^m$, assume that \\
$\int_{\mathbb R}\exp\left(-Q_t^{m, \pi}(v;\xi,w)/\lambda\right)\, \mathrm dv <\infty$.
Define the improved Gibbs policy $\pi'=\{\pi_t'\}_{t=0}^{T-1}$ by
\begin{align}\label{definpitrho}
\pi_t'(u; \xi)
:=
\frac{
\exp\left(-Q_t^{m, \pi}(u;\xi,w)/\lambda\right)
}{
\displaystyle \int_{\mathbb R}\exp \left(-Q_t^{m, \pi}(v;\xi, w)/\lambda\right)\, \mathrm dv
}.
\end{align}
If the normalizing integrals are positive and the resulting policy $\pi'$ is admissible, then $J_t^{m,\pi'}(\xi;w)\le J_t^{m,\pi}(\xi;w)$.
\end{theorem}

The convergence of the resulting policy iteration is stated below, with the proof provided in \ref{thm:policy_convergence00}.

\begin{theorem}[Policy convergence]\label{thm:policy_convergence}
Let $\pi^{(0)}$ be any admissible initial feedback control policy. Fix $w$ and suppose that the evaluations are finite and each Gibbs update has a positive finite normalizing integral and defines an admissible feedback policy. For each $n\ge0$, define recursively:
\begin{enumerate}[label=(\roman*)]
\item policy evaluation
\begin{align*}
& J_T^{m,\pi^{(n)}}(\xi;w)=(x_T^{\pi^{(n)}}-l_T-w)^2-(w-d)^2, \qquad  \quad t=T,\\
& J_t^{m,\pi^{(n)}}(\xi;w)
=
(\mathcal T_t^{m,\pi^{(n)}}J_{t+1}^{m,\pi^{(n)}})(\xi),
\qquad \qquad \qquad \ \
t=T-1,T-2,\cdots,0;
\end{align*}
\item policy improvement
\begin{align}\label{eq:soft_policy_iteration_corrected}
\pi_t^{(n+1)}(u; \xi)
=
\frac{
\exp\big(-Q_t^{m,\pi^{(n)}}(u;\xi,w)/\lambda\big)
}{
\displaystyle
\int_{\mathbb R}
\exp\big(-Q_t^{m,\pi^{(n)}}(v;\xi,w)/\lambda\big)\,\mathrm dv
},
\qquad t=0,1,\cdots,T-1.
\end{align}
\end{enumerate}
Then $J_t^{m,\pi^{(n+1)}}(\xi;w)\le J_t^{m,\pi^{(n)}}(\xi;w)$ for $n\ge 0$. For every $t<T$, $J_t^{m,\pi^{(n)}}(\xi;w)=J_t^{m,*}(\xi;w)$ for $n\ge T-t$.

Moreover, for $n\ge T-t$, the density $\pi_t^{(n)}(\cdot;\xi)$ agrees almost everywhere with $\pi_t^{m,*}(\cdot;\xi)$, where
\begin{align}\label{pitstaruxi}
\pi_t^{m,*}(u; \xi)
=
\frac{
\exp\big(-Q_t^{m,*}(u;\xi,w)/\lambda\big)
}{
\displaystyle
\int_{\mathbb R}
\exp\big(-Q_t^{m,*}(v;\xi,w)/\lambda\big)\,\mathrm dv
},
\end{align}
with $Q_t^{m,*}(u;\xi,w)
:=
\mathbb E\left[
J_{t+1}^{m,*}(\xi_{t+1}^m;w)
\;\middle|\;
\xi_t^m=\xi,\ u_t=u
\right]$.
\end{theorem}

\begin{remark}\label{rem:gaussian_special_case_corrected}
For a given case $m$ and iteration $n$, whenever the mapping $u\mapsto Q_t^{m,\pi^{(n)}}(u;\xi,w)$ is quadratic with positive leading coefficient, the general Gibbs policy-improvement step \eqref{definpitrho} reduces to a Gaussian update, that is, to a recursive update of the mean and variance parameters.
Accordingly, the explicit Gaussian formulas obtained in the complete-information case, in the partially observed case under quadratic closure, and in the no-learning benchmark, cf. \eqref{pi32a}, \eqref{oppoalp}, and \eqref{pi316}, respectively, should be viewed as particular instances of \eqref{definpitrho}.
\end{remark}

\subsection{RL Algorithm Design}
Motivated by Gibbs policy improvement, we implement a finite-horizon entropy-regularized actor-critic approximation with stochastic gradient updates (SGD). The exact policy-iteration result applies to the Bellman models above. In the simulation, the learned no-learning policy is trained on the original switching paths with the prior as its regime input. The algorithm involves five steps:
\emph{initialization}, \emph{sampling}, \emph{policy evaluation}, \emph{policy gradient}, and \emph{self-correcting scheme for learning the Lagrange multiplier $w$}.

\emph{Step 1: Initialization.}
In view of the theoretical results in Theorems \ref{thm32a}-\ref{negtilpthm}, a quadratic structure is adopted for the value function and a Gaussian parameterization for the policy. Accordingly, for each $m\in\{c,p,s\}$, define the parameterized value function
\begin{align}\label{eq:alg_value_param_structured}
J_t^{m,\eta}(x,l,z;w)
&=
a_t^{m,\eta}(z)\,x^2
+\big(b_{t,0}^{m,\eta}(z) + b_{t,1}^{m,\eta}(z) w + b_{t,2}^{m,\eta}(z) l\big) \,x
 \nonumber
\\
& \ \
+ \big(c_{t,0}^{m,\eta}(z)+ c_{t,1}^{m,\eta}(z) w + c_{t,2}^{m,\eta}(z) l \big)\, l
+ e_{t,0}^{m,\eta}(z) + e_{t,1}^{m,\eta}(z)\, w + e_{t,2}^{m,\eta}(z) \,w^2,
\end{align}
with terminal condition $J_T^{m,\eta}(x, l, z;w) = (x-l-w)^2-(w-d)^2$, where the coefficient functions are parameterized by low-order polynomial bases in the state variable $z$ and the remaining horizon $T-t$, more precisely, for $\ell=0,1,2$, denote $\eta:=(\eta_a, \eta_{b, \ell}, \eta_{c, \ell}, \eta_{e, \ell})$ with $\eta_a:=(\eta_{a, r, q})_{r, q}$, $\eta_{\cdot, \ell}:=(\eta_{\cdot, \ell, r, q})_{r, q}$, and
\begin{align*}
& a_{t}^{m,\eta}(z)
=
\exp\big\{\sum_{r=0}^{\mathds R_{a}}\sum_{q=0}^{\mathds Q_{a}}
\eta_{a,r,q} z^r (T-t)^q\big\},
\quad
b_{t,\ell}^{m,\eta}(z)
=
-\exp\big\{\sum_{r=0}^{\mathds R_{b,\ell}}\sum_{q=0}^{\mathds Q_{b,\ell}}
\eta_{b,\ell, r,q} z^r (T-t)^q\big\}, \\
& c_{t,\ell}^{m,\eta}(z)
= \sum_{r=0}^{\mathds R_{c,\ell}}\sum_{q=0}^{\mathds Q_{c,\ell}}
\eta_{c,\ell, r,q} z^r (T-t)^q,
\qquad \quad \ \
e_{t,\ell}^{m,\eta}(z)
= \sum_{r=0}^{\mathds R_{e,\ell}}\sum_{q=0}^{\mathds Q_{e,\ell}}
\eta_{e,\ell, r,q} z^r (T-t)^q.
\end{align*}
In the complete-information case, since $z=\varepsilon_t\in\{1,2\}$, the above polynomial dependence can equivalently be replaced by regime-specific coefficients.
We then parameterize the actor by a Gaussian density
\begin{align}\label{eq:alg_policy_param}
\pi_{t}^{m,\phi}(u; x,l,z,w)
= \mathcal N\!\left(
u\;\middle|\;
\mu_t^{m,\phi}(x,l,z,w),\,
\sigma_t^{m, \phi}(x,l,z,w)
\right),
\end{align}
where $\phi:=(\phi_{\alpha, \ell}, \phi_{\beta})$ with $\phi_{\alpha, \ell}:=(\phi_{\alpha, \ell, r,q})_{r, q}$ and $\phi_{\beta}:=(\phi_{\beta, r, q})_{r, q}$, $\ell = 0,1,2$, \\ $\alpha_{t,\ell}^{m,\phi}(z)
= \sum_{r=0}^{\mathds R_{\alpha,\ell}}\sum_{q=0}^{\mathds Q_{\alpha,\ell}}
\phi_{\alpha,\ell, r,q} z^r (T-t)^q$,
$\beta_{t}^{m,\phi}(z)
= \sum_{r=0}^{\mathds R_{\beta}}\sum_{q=0}^{\mathds Q_{\beta}}
\phi_{\beta, r,q} z^r (T-t)^q$, and
\begin{align}
&\mu_t^{m, \phi}(x,l,z,w)
= b_{t,0}^{m,\eta}(z) \,
\alpha_{t,0}^{m,\phi}(z)\,x
-\big(\frac{b_{t,1}^{m,\eta}(z)}{a_t^{m,\eta}(z)} \,
e^{\alpha_{t,1}^{m,\phi}(z)}\, w+ \frac{b_{t,2}^{m,\eta}(z)}{a_t^{m,\eta}(z)} \,
e^{\alpha_{t,2}^{m,\phi}(z)}\, l), \label{eq:alg_policy_mean} \\
&
\sigma_t^{m, \phi}(x,l,z,w)
= \frac{\lambda}{2 a_t^{m,\eta}(z)} e^{\beta_{t}^{m,\phi}(z)}. \label{eq:alg_policy_var}
\end{align}
For simplicity, set $\mathds R_a =\mathds R_{b, \ell}= \mathds R_{c, \ell}= \mathds R_{e, \ell}=\mathds R_{\alpha, \ell}=\mathds R_{\beta}:= \mathds R$ and $\mathds Q_{a}= \mathds Q_{b, \ell} = \mathds Q_{c, \ell} = \mathds Q_{e, \ell}=\mathds Q_{\alpha, \ell}=\mathds Q_{\beta}:= \mathds Q$. Again, the coefficient functions in \eqref{eq:alg_policy_mean}-\eqref{eq:alg_policy_var} are parameterized by low-order polynomial bases in $z$ and $T-t$. The exponential parameterization is used to guarantee the positivity of the variance during training. The entropy of the Gaussian policy \eqref{eq:alg_policy_param} admits the closed-form expression
\begin{align}\label{eq:alg_entropy}
\mathcal H(\pi_{t}^{m, \phi})
&:=
-\int_{\mathbb R}\pi_{t}^{m,\phi}(u; x,l,z,w)\ln \pi_{t}^{m, \phi}(u; x,l,z,w)\,\mathrm{d}u
=
\frac12\ln \bigl(2\pi e\,\sigma_t^{m,\phi}(x,l,z,w)\bigr) \nonumber\\
&= -\frac{1}{2}\ln \bigl(\frac{a_t^{m,\eta}(z)}{\pi \lambda}\bigl) + \frac{1}{2}(\beta_{t}^{m,\phi}(z)+1).
\end{align}
The coefficient functions $a_{t}^{m,\eta}(z)$, $b_{t,\ell}^{m,\eta}(z)$, $c_{t,\ell}^{m,\eta}(z)$, and $e_{t, \ell}^{m,\eta}(z)$ are determined by the value-function parameter vection $\eta$, whereas $\alpha_{t,\ell}^{m,\phi}(z)$ and $\beta_t^{m, \phi}(z)$ are parameterized by the policy parameter vector $\phi$. The vectors $\eta$, $\phi$, together with the Lagrange multiplier $w$, are learned during training and initialized arbitrarily.
The parameterizations \eqref{eq:alg_value_param_structured} and \eqref{eq:alg_policy_param} retain the quadratic and Gaussian structure of the complete-information and auxiliary prior-mixture solutions, derived in Theorems \ref{thm32a} and \ref{negtilpthm}. The finite polynomial basis approximates the coefficient functions. Under partial information, the Gaussian actor provides a structured approximation to the Gibbs policy in Theorem \ref{them1122}.

\emph{Step 2: Sampling.}
At each training iteration, we generate a single sample trajectory (or a batch of $K$ trajectories and use the resulting sample average) under the current policy, which consists of $T$ time steps. More precisely, at time $t$, given the current hidden regime $\varepsilon_t=i$, current wealth-liability state $(x_t^\pi,l_t)$, and action $u_t$, the next observation $Y_{t+1}:=(x_{t+1}^\pi,l_{t+1})$ is generated from the conditional law $Y_{t+1}\sim G_t^{(\varepsilon_t)}(\cdot\,;x_t^\pi,l_t,u_t)$.
\begin{enumerate}[label=(\roman*)]
\item In the complete-information case, sample $u_t\sim\pi_{t}^{c, \phi}(\cdot; \xi_t^c,w)$, generate the next observation $Y_{t+1}:=(x_{t+1}^{\pi}, l_{t+1})$, sample $\varepsilon_{t+1}\sim P(\varepsilon_t,\cdot)$, and simulate $\xi_{t+1}^c=(x_{t+1}^\pi,l_{t+1},\varepsilon_{t+1})$.
\item In the partially observed case, sample $u_t\sim\pi_{t}^{p, \phi}(\cdot; \xi_t^p,w)$, generate the next observation $Y_{t+1}:=(x_{t+1}^{\pi}, l_{t+1})$, update the posterior belief $\hat p_{t+1}$ by the Bayes recursion \eqref{titer}, sample $\varepsilon_{t+1}\sim P(\varepsilon_t,\cdot)$, and simulate $\xi_{t+1}^p=(x_{t+1}^\pi,l_{t+1},\hat p_{t+1})$.
\item In the no-learning benchmark, sample $u_t\sim\pi_{t}^{s, \phi}(\cdot; \xi_t^s,w)$, generate the next observation $Y_{t+1}:=(x_{t+1}^{\pi}, l_{t+1})$, update $\tilde p_{t+1}$ by \eqref{tildep00},  sample $\varepsilon_{t+1}\sim P(\varepsilon_t,\cdot)$, and simulate $\xi_{t+1}^s=(x_{t+1}^\pi,l_{t+1},\tilde p_{t+1})$.
\end{enumerate}
Accordingly, for each $m\in\{c,p,s\}$, one sampled episode is
$
\mathcal D^{m}
=
\bigl\{(\xi_t^m,u_t,\xi_{t+1}^m)\bigr\}_{t=0}^{T-1} \ \ \mathrm{with} \ \ \xi_t^m=(x_t^\pi,l_t,z_t^m)$. 

\emph{Step 3: Policy evaluation (PE).}
For a Markov policy $\pi$ evaluated in its Bellman model, the following process is a martingale with respect to the pre-decision history of states and actions,
$M_t^m:={J}_{t}^{m, \pi}(\xi; w)-\lambda\sum_{k=0}^{t-1}\mathcal{H}(\pi_k)$. 
To fit the critic parameters $\eta$ in \eqref{eq:alg_value_param_structured} for a fixed actor $\pi_t^{m,\phi}$, we consider 
$M_t^{m, \eta, \phi}:={J}^{m, \eta}_t(\xi; w)
-\lambda\sum_{k=0}^{t-1}\mathcal{H}(\pi_k^{m, \phi})$. 
Assume that the terminal cost, fitted critic values at each date, and accumulated entropy are square-integrable. For policy evaluation we minimize the martingale loss (ML) function $\mathrm{ML}^m(\eta;\phi,w)$
\begin{align*}
\mathrm{ML}^m(\eta; \phi, w)
&:=\frac{1}{2} \mathbb{E}
\big[\sum_{t=0}^{T-1} | M_T^{m, \eta, \phi} -M_t^{m, \eta, \phi}| ^2 \Delta t \big]\\
&=\frac{1}{2} \mathbb{E}
\big[\sum_{t=0}^{T-1} \big( J_{T}^{m, \eta}(\xi_T^m; w)
-J_{t}^{m, \eta}(\xi_t^m; w)
-\lambda\sum_{k=t}^{T-1}\mathcal{H}(\pi_k^{m, \phi})\big)^2 \Delta t \big],
\end{align*}
where
$J_{T}$ is given by \eqref{term098}. After sampling, the above expectation is approximated by the empirical batch $\mathcal D^m$ (or average over $K$ batches), and utilize SGD algorithm to design the parameter updating rule for $\eta$. Specifically, for $(\xi_{t}^{m},{u}_{t}, \xi_{t+1}^{m})\in \mathcal D^m$ with $\xi_t^m = (x_t^{\pi^{m, \phi}}, l_t, z_t^m)$, the gradient of $\mbox{ML}^m$ with respect to $\eta$ is
\begin{align}\label{eq:grad_ML_unified}
\nabla_\eta \mathrm{ML}^m(\eta;\phi,w)
&=-\mathbb{E}
\Big[\sum_{t=0}^{T-1} \Big( \big( J_{T}^{m, \eta}(\xi_T^m; w)
-J_{t}^{m,\eta}(\xi_t^m; w)
-\lambda\sum_{k=t}^{T-1}\mathcal{H}(\pi_k^{m, \phi})\big) \Delta t
\cdot \nabla_\eta J_t^{m,\eta}(\xi_t^{m};w) \Big) \Big],
\end{align}
where the sampled trajectories, executed policy, policy entropy, and $w$ are held fixed during this step, $\nabla_\eta\mathcal H(\pi_k^{m,\phi})=0$ during policy evaluation, and for $r=0,1,\cdots, \mathds R$, $q=0,1,\cdots, \mathds Q$,
\begin{align*}
& \frac{\partial J_{t}^{m, \eta}(\xi; w)}{\partial \eta_{a,r,q}}= a_t^{m, \eta}(z) x^2 z^r (T-t)^q, \\
& \frac{\partial J_{t}^{m, \eta}(\xi; w)}{\partial \eta_{b,\ell, r,q}}
= \big(b_{t,0}^{m,\eta}(z)\mathbf{1}_{\{\ell=0\}} + b_{t,1}^{m,\eta}(z) w \mathbf{1}_{\{\ell=1\}} + b_{t,2}^{m,\eta}(z) l \mathbf{1}_{\{\ell=2\}} \big) \,x z^r (T-t)^q, \\
& \frac{\partial J_{t}^{m, \eta}(\xi; w)}{\partial \eta_{c,\ell, r,q}}
= \big(\mathbf{1}_{\{\ell=0\}} + w \mathbf{1}_{\{\ell=1\}} +  l \mathbf{1}_{\{\ell=2\}} \big)\, l z^r (T-t)^q,\\
& \frac{\partial J_{t}^{m, \eta}(\xi; w)}{\partial \eta_{e,\ell, r,q}}
= \big(\mathbf{1}_{\{\ell=0\}} + w \mathbf{1}_{\{\ell=1\}} +  w^2 \mathbf{1}_{\{\ell=2\}} \big)\, z^r (T-t)^q, \\
& \frac{\partial \mathcal{H}(\pi_k^{m, \phi})}{\partial \eta_{a,r,q}}
= - \frac{1}{2} z^r (T-t)^q,
\qquad
\frac{\partial \mathcal{H}(\pi_k^{m, \phi})}{\partial \eta_{b, \ell,r,q}}
=  \frac{\partial \mathcal{H}(\pi_k^{m, \phi})}{\partial \eta_{c, \ell,r,q}}
=\frac{\partial \mathcal{H}(\pi_k^{m, \phi})}{\partial \eta_{e, \ell,r,q}} =0.
\end{align*}
The value-function parameter $\eta$ is then updated by SGD, $\eta \leftarrow \eta-\alpha_{\eta}^m \,\nabla_\eta \mathrm{ML}^m(\eta;\phi,w)$,
with learning rate $\alpha_{\eta}^m>0$. The implementation attempts one critic update per training iteration, with backtracking on the batch loss.

\emph{Step 4: Policy gradient.}
The policy \eqref{eq:alg_policy_param} is updated using the gradient of the objective function, $G^m(\phi; w)$, as defined by Theorem 5 in \cite{Jia2022},
\begin{align*}
G^m(\phi; w):=\mathbb{E}&\Big[\sum_{t=0}^{T-1}  \Big(\nabla_{\phi} \ln \pi_t^{m, \phi}(u; \xi, w)
\big(J_{t+1}^{\pi^{m, \phi}}(\xi_{t+1}^m; w)
- J_{t}^{\pi^{m, \phi}}(\xi_t^m; w) -\lambda \mathcal{H}(\pi_t^{m, \phi}) \big)
-\lambda \nabla_{\phi} \mathcal{H}(\pi_t^{m, \phi}) \Big) \Big],
\end{align*}
where $u$ is the sampled action from $\pi_t^{m, \phi}$. By replacing $J_{t}^{\pi^{m, \phi}}$ by its parameterized function ${J}^{m,\eta}_t$, we obtain
\begin{align}\label{phi1}
G^m(\phi; \eta, w)
:=\mathbb{E}&\Big[\sum_{t=0}^{T-1} \Big(\nabla_{\phi} \ln \pi_t^{m, \phi}(u; \xi, w)
\big(J_{t+1}^{m, \eta}(\xi_{t+1}^m; w)
- J_{t}^{m, \eta}(\xi_t^m; w)
-\lambda \mathcal{H}(\pi_t^{m, \phi})\big)-\lambda \nabla_{\phi} \mathcal{H}(\pi_t^{m, \phi}) \Big) \Big],
\end{align}
where
\begin{align*}
&\frac{\partial \ln \pi_t^{m, \phi}(u; \xi,w)}{\partial \phi_{\alpha,0, r, q}}
=\frac{1}{\sigma_t^{m, \phi}(x,l,z,w)}
\big(u - \mu_t^{m, \phi}(x,l,z,w)\big) b_{t,0}^{m,\eta}(z) x \cdot z^r (T-t)^q, \\
&\frac{\partial \ln \pi_t^{m, \phi}(u; \xi,w)}{\partial \phi_{\alpha,1, r, q}}
= - \frac{1}{\sigma_t^{m, \phi}(x,l,z,w)}
\big(u - \mu_t^{m, \phi}(x,l,z,w)\big) \frac{b_{t,1}^{m,\eta}(z)}{a_t^{m,\eta}(z)} \,
e^{\alpha_{t,1}^{m,\phi}(z)}\, w \cdot z^r (T-t)^q,\\
&\frac{\partial \ln \pi_t^{m, \phi}(u; \xi, w)}{\partial \phi_{\alpha,2, r, q}}
= - \frac{1}{\sigma_t^{m, \phi}(x,l,z,w)}
\big(u - \mu_t^{m, \phi}(x,l,z,w)\big) \frac{b_{t,2}^{m,\eta}(z)}{a_t^{m,\eta}(z)} \,
e^{\alpha_{t,2}^{m,\phi}(z)}\, l \cdot z^r (T-t)^q,\\
&\frac{\partial \ln \pi_t^{m, \phi}(u; \xi, w)}{\partial \phi_{\beta, r, q}}
=\Big(\frac{a_t^{m,\eta}(z)}{\lambda}\big(u - \mu_t^{m, \phi}(x,l,z,w)\big)^2 e^{-\beta_{t}^{m,\phi}(z)}
- \frac{1}{2}\big) \cdot z^r (T-t)^q, \\
&\frac{\partial \mathcal{H}(\pi_t^{m, \phi})}{\partial \phi_{\alpha, 0, r, q}}=
\frac{\partial \mathcal{H}(\pi_t^{m, \phi})}{\partial \phi_{\alpha, 1, r, q}}=
\frac{\partial \mathcal{H}(\pi_t^{m, \phi})}{\partial \phi_{\alpha, 2, r, q}}=0,\ \
\frac{\partial \mathcal{H}(\pi_t^{m, \phi})}{\partial \phi_{\beta, r, q}}=\frac{1}{2} z^{r} (T-t)^{q}.
\end{align*}
The policy-parameter update is $\phi \leftarrow
\phi-\alpha_{\phi}^m \mathrm{G}^m(\phi;\eta,w)$,
where $\alpha_{\phi}^m>0$ is the policy (actor) learning rate.

\emph{Step 5: Self-correcting scheme for learning the Lagrange multiplier $w$.}
Update $w$ by $w \leftarrow w - \alpha_w^m ({x}_T^{\pi} - {l}_T -d)$, where $\alpha_w^m>0$ is the learning rate.
This scheme works as follows: if the empirical terminal surplus falls below the target $d$, the self-correcting rule raises $w$ to increase the mean of our strategy \eqref{eq:alg_policy_param}, guiding the terminal surplus to increase to the expected level. Generally, $w$ should be updated gradually to minimize variability in its learning process. We implement the update for $w$ after every $N$ training iterations.

In the implementation, let $\bar S_k$ denote the average terminal surplus of the training batch at iteration $k$. At a multiplier update, the empirical surplus estimate is the average of the most recent $N$ batch means.
The algorithms corresponding to the complete-information case, the partially observed case, and the no-learning benchmark are hereafter denoted by CoEMV, PoEMV-1, and PoEMV-2, respectively.
\begin{algorithm}[H]\label{Algorithm1}
        \caption{Multi-period CoEMV/PoEMV-1/PoEMV-2 Algorithm}
        \KwIn{Market, learning rates $\alpha_{\eta}^m$, $\alpha_{\phi}^m$, $\alpha_{w}^m$, initial wealth $x_0$, initial liability $l_0$, initial market state $\varepsilon_0$, initial distribution of the market state $\hat p_0$ and $\tilde p_0$, target payoff $d$, investment horizon $T$, time step $\Delta t$, exploration weight $\lambda$, transition probability matrix of the market state $P=(P_{ij})_{2\times 2}$, number of iterations $n$, batch size $B$ and sample average size $N$.}
        \KwOut{Learned portfolio strategy and approximate value function}
        initialization $\eta$, $\phi$ and $w$\\
        \For{$k=1$ to $n$}{
             \For{$t=0$ to $T-1$}
             {Sample $B$ episodes under the current policy $\pi^{m,\phi}$, each consisting of the transitions $(\xi_t^m,u_t^\pi,\xi_{t+1}^m)$ for $t=0,\ldots,T-1$. Record their mean terminal surplus $\bar S_k$.}
        \textbf{Policy Evaluation Update}\\
        {
        Update $\eta \leftarrow \eta-\alpha_{\eta}^m\,\nabla_\eta \mathrm{ML}^m(\eta;\phi,w)$ using \eqref{eq:grad_ML_unified}}\\
         \textbf{Policy Gradient Update}\\
        {Update $\phi \leftarrow \phi-\alpha_{\phi}^m \mathrm{G}^m(\phi;\eta,w)$ using \eqref{phi1}}\\
        Update $\pi^{m, \phi}\leftarrow \pi^{m, \phi}(u; x, l, z, w)$ that given by \eqref{eq:alg_policy_param}\\
        \textbf{Lagrange multiplier Update}\\
            \eIf{ k mod N=0}
            {
                    Update $w\leftarrow w-\alpha_w^m \big(\frac{1}{N}\sum_{j=k-N+1}^k \bar{S}_j-d\big)$ \\
            }
            {
                $w\leftarrow w$ \\
            }
    }
\end{algorithm}

\section{Numerical analysis}\label{sec5}

Direct evaluation of the Bellman recursion in Theorem \ref{them1122} requires the transition and observation laws. Based on the derived feedback structure and the multi-period approach of \cite{Cui2023}, we use PoEMV-1 to learn a parameterized approximation. The filter uses the specified simulation parameters or market parameters estimated from the historical training sample. Additionally, we compare the learned PoEMV-1 policy with policies that use the prior (PoEMV-2) or rely on complete information (CoEMV).

\subsection{Simulation study}\label{simu}
We first specify several parameter values. The investment horizon is set to $T=10$ years with daily time step $\Delta t=1/252$ (cited in \cite{Wu2024}), giving $T=2520$ decision periods. The initial wealth and liability are $x_0=1$ and $l_0=0.1$, respectively, and the terminal target is $d=8$ (cited in \cite{Wang2020}). The initial probability of regime 1 is $\hat p_0=0.3$, and the daily transition matrix of the hidden market regime is specified by $P_{11}=0.9986$, $P_{12}=0.0014$, $P_{21}=0.0114$ and $P_{22}=0.9886$. For the risk-free component, the annual gross returns are $e^0(1)=1.20$ and $e^0(2)=1.05$, which are taken from \cite{Liang2015} and reflect the fact that asset returns are higher in a bull market, and converted into one-step gross returns through $e^0_{\Delta t}(i)=(e^0(i))^{\Delta t}$.
Since the skewed $t$ distribution is proved by \cite{Theodossiou1998} to effectively capture skewness and tail risk, as well as to provide an excellent fit to the empirical distribution of financial data, we use this distribution for simulating the net risky-asset return $e_t^1(\varepsilon_t)-1$ with parameters that vary by $\varepsilon_t$: if $\varepsilon_t = 1$, then the parameters borrowed from \cite{Cui2023} are annualized mean $0.5$, volatility $0.12$, 10 degrees of freedom and skewness $0.1$; if $\varepsilon_t=2$, then the parameters shift to a return of $0.06$, volatility $0.18$, with the same degrees of freedom and skewness.
The one-step mean and volatility are obtained by multiplying their annualized counterparts by $\Delta t$ and $\sqrt{\Delta t}$, respectively. 
The one-step net liability returns $q_t(1)-1$ and $q_t(2)-1$ follow normal distributions $\mathcal N(0.05\Delta t,0.036\Delta t)$ and $\mathcal N(0.01\Delta t,0.072\Delta t)$, respectively.

\begin{figure}[htbp]
	\centering
	\includegraphics[width=5in]{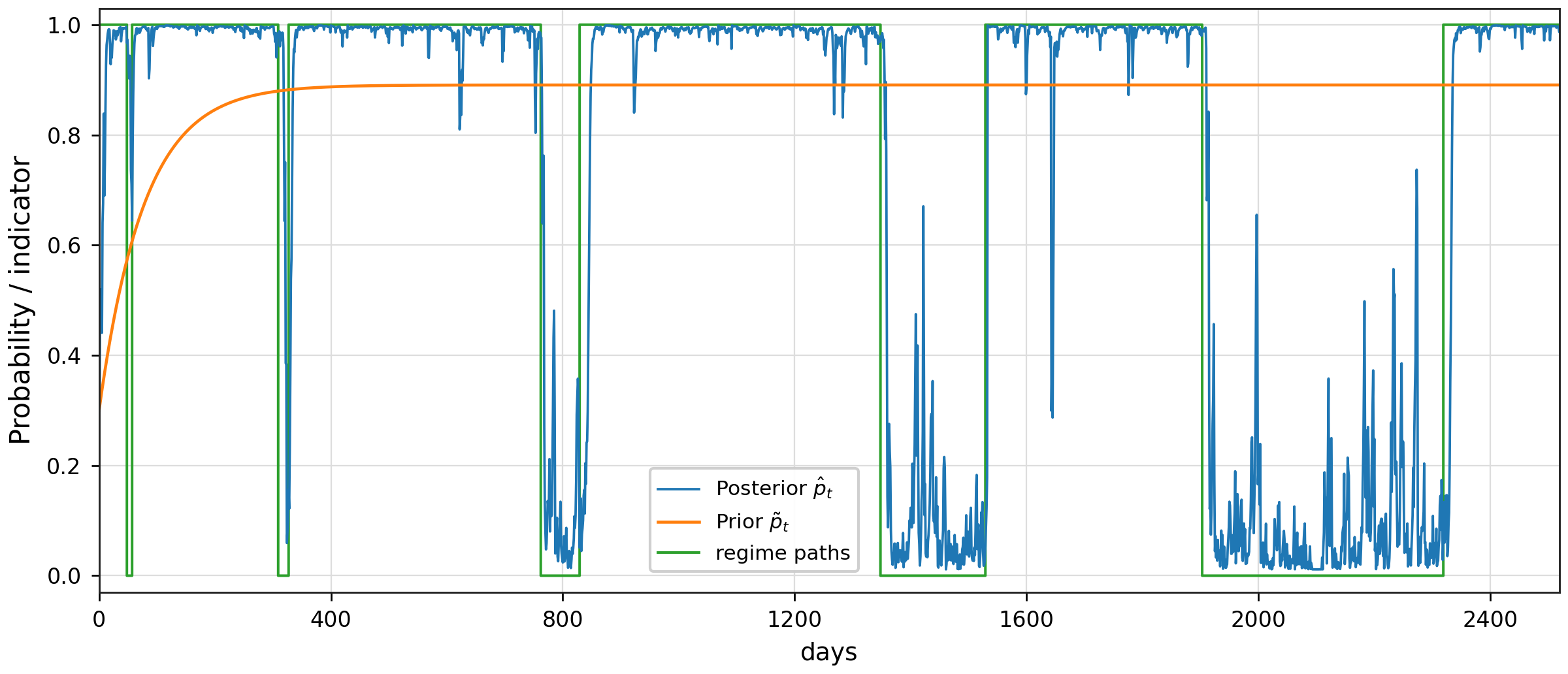}
	\caption{Evaluated market regimes over 10 years.}
	\label{fig:1}
\end{figure}

The market regime is evaluated in three ways. First, in the partially observed case, the posterior belief $\hat p_t$ is obtained by the Bayesian update and prediction in Lemma \ref{p2t}, implemented in log-odds form conditional on the executed action and the state observation $Y_{t+1}=(x_{t+1},l_{t+1})$. In this recursion, the likelihood under each candidate regime is computed from the risky-asset and liability returns implied by the observed state transition. Second, the no-learning benchmark uses the prior path $\tilde p_t$ in \eqref{tildep00}, which is propagated only by the transition matrix and is not updated by the realized observation. Third, the complete-information benchmark uses the realized regime path directly. Fig.\ref{fig:1} reports, over a 10-year horizon, a posterior path $\hat p_t$, the corresponding no-learning prior path $\tilde p_t$, and the realized regime process.

We next solve the entropy-regularized multi-period EMV problem \eqref{min22} by the finite-horizon actor-critic scheme described in Algorithm \ref{Algorithm1}. In all three information structures, namely PoEMV-1, PoEMV-2 and CoEMV, the value function is parameterized by the quadratic form in \((x,l,w)\), with coefficient functions expanded in the basis \(z^r[(T-t)/T]^q\), where \(r,q=0,1,2\) and $\mathds R=\mathds Q=2$.
The Gaussian actor retains the parametric structure derived from the entropy-regularized policy-improvement step, and the annual entropy intensity is set to \(2\), giving the one-step weight \(\lambda=2\Delta t\). 
The learning rates for the value-function parameters are
\(
(\alpha_{\eta_a},\alpha_{\eta_b},\alpha_{\eta_c},\alpha_{\eta_e})
=
(7.5\times 10^{-4},\,7.5\times 10^{-4},\,5.0\times 10^{-4},\,5.0\times 10^{-4}),
\)
and those for the policy parameters are
\(
(\alpha_{\phi_\alpha},\alpha_{\phi_\beta})
= (2.0\times 10^{-5},\,2.0\times 10^{-6})
\).
The learning rate for the Lagrange multiplier is \(\alpha_w=4\times 10^{-3}\).

The numerical implementation uses a terminal-consistent, state-feedback warm start for the value-function and policy parameters to enhance numerical stability over the long investment horizon. For PoEMV-1, candidate initial policies are obtained by blending the common initial feedback with the posterior-weighted reference feedback and projecting the resulting feedback onto the actor basis. Among the candidates calibrated to the target mean \(d=8\) on \(128\) fitting paths, we retain the one with the lowest variance of terminal net wealth. Given the initial policy for each algorithm, we then calibrate \(w_0\) on \(64\) simulated paths to give a mean terminal net wealth of \(6\). Under the simulation setting used here, this gives \(w_0\approx5.9388\) for PoEMV-1 and \(w_0\approx6.2113\) for both PoEMV-2 and CoEMV.

In each training iteration, a batch of \(4\) independent 10-year episodes is generated under the current Gaussian policy, with the Gaussian variance constrained to \([10^{-8},0.02]\) and the sampled action clipped to \([-20,20]\) before execution. The critic minimizes the martingale loss by gradient descent with the policy and its entropy held fixed. After the first \(750\) critic-only iterations, the actor update follows the policy-gradient direction in \eqref{phi1}, using a slowly updated target critic and a batchwise leave-one-trajectory-out baseline to reduce variance, with a gradual increase in the nominal step size, subsequent power-law decay, and KL-controlled backtracking for numerical stability. The critic and actor updates use the averages of their respective gradient estimates over the current batch. Following the self-correcting rule in \cite{Wu2024}, the multiplier is updated every \(N=10\) iterations after the critic-only phase, using the mean terminal net wealth over the most recent \(10\) batches. All three algorithms are run for \(20{,}000\) iterations (resulting in \(80{,}000\) episodes each).

\begin{figure}[htbp]
	\centering
	\includegraphics[width=6in]{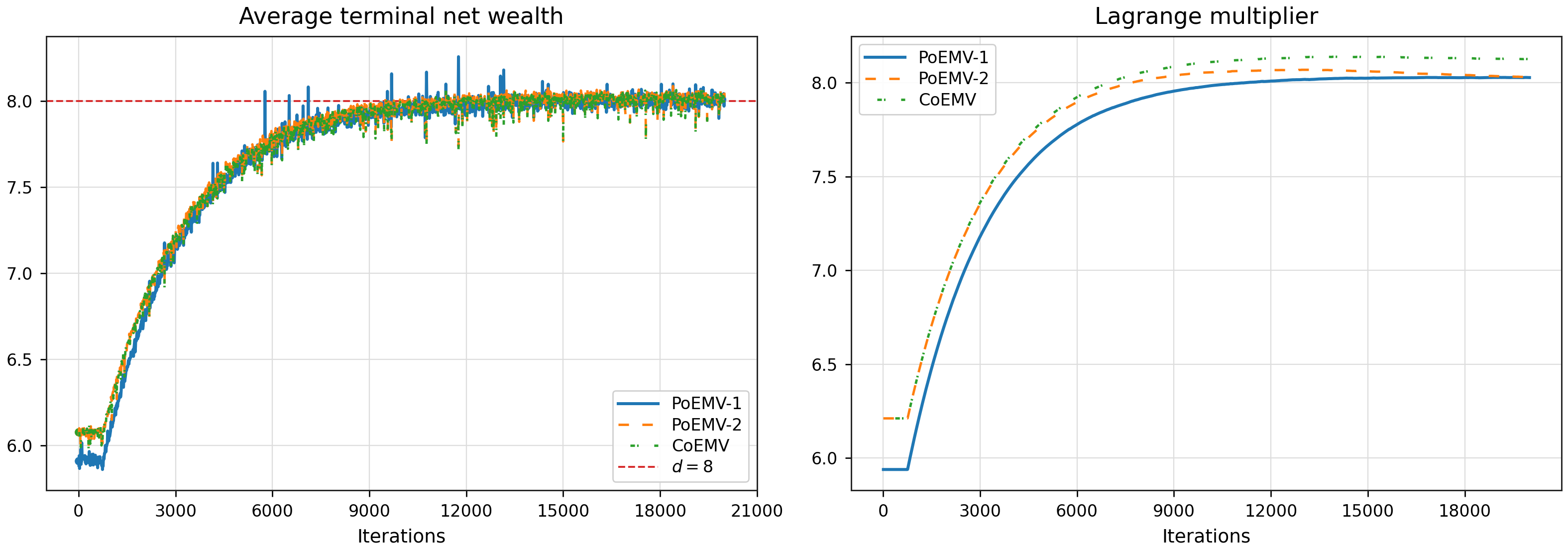} \ \
	\caption{The final performance (averaged over every 10 iterations) in the training process.}
	\label{fig:2}
\end{figure}

The training performance is summarized by block averages over every \(10\) iterations. Fig.\ref{fig:2} displays the block-average terminal net wealth \(x_T-l_T\) and the block-average Lagrange multiplier \(w\). The block averages for all three algorithms and the out-of-sample results reported below use the seed \(34567\).

We further assess the out-of-sample performance of the EMV strategies using three measures: the sample terminal mean net wealth, the sample variance of terminal net wealth across testing paths, and the 10-year Sharpe ratio, defined as
\(
(\mathrm{Mean}-x_0)/\sqrt{\mathrm{Variance}},
\)
where \(x_0=1\). We generate \(1000\) independent 10-year testing paths from the same regime-switching asset-liability model \eqref{st4}, which are not used during training and validation, and apply all policies to the same paths.
For out-of-sample evaluation, at each step we draw an action from the fitted Gaussian policy and clip it to \([-20,20]\), with the Gaussian variance constrained to \([10^{-8},0.02]\). The initial policy is evaluated at iteration \(0\).
Every \(1000\) training iterations thereafter, the current and averaged policies are evaluated on a fixed set of \(256\) independent validation paths. Among the candidates whose validation mean is within \(0.10\) of \(d\), we select the one with the lowest validation variance. If none satisfies this tolerance, we select the candidate with the smallest absolute mean-target gap.
The selected policy and its associated multiplier are then used for testing. For the reported seed, the averaged PoEMV-1 policy is selected at iteration \(6000\), and the current PoEMV-2 and CoEMV policies are selected at iteration \(20000\). All three learned policies use recalibrated multipliers.

For comparison, we also evaluate the theoretical optimal benchmark policies associated with the three information structures. In the partially observed case, the benchmark denoted by Optimal-1 is obtained from a posterior-projected quadratic Gibbs policy, implemented numerically with the posterior belief recursively updated by the action-dependent state-observation filter.
In the no-learning and complete-information cases, the benchmarks denoted by Optimal-2 and Optimal-3 are generated from the explicit feedback controls $\pi_t^{sub*}$ in \eqref{pi316} and ${\pi}_{t}^{c*}$ in \eqref{pi32a}, respectively. For these theoretical benchmarks, the Lagrange multiplier \(w\) is calibrated separately by a one-dimensional root-finding procedure so that the simulated mean terminal net wealth matches the target level \(d=8\).
The resulting out-of-sample statistics are reported in Table \ref{Table 1}, where the mean and variance are computed from the terminal net wealth \(x_T-l_T\) over the \(1000\) testing paths, and the Sharpe ratio is calculated from these two sample quantities.

\begin{table}[!ht]
\caption{Out-of-sample performance.}
\label{Table 1}
\centering
\begin{tabular}{l| c c c}
\toprule
$\ $ & $\text{Mean}$ & $\text{Variance}$ & $\text{Sharpe ratio}$ \\
\cline{1-4}
$\text{PoEMV-1}$ & $7.9954$ & $0.0322$ & $38.9836$\\
$\text{Optimal-1}$ & $8.0104$ & $0.0207$ & $48.6538$\\
\midrule
$\text{PoEMV-2}$ & $8.0074$ & $0.0443$ & $33.3068$\\
$\text{Optimal-2}$ & $7.9900$ & $0.1172$ & $20.4191$\\
\midrule
$\text{CoEMV}$ & $8.0064$ & $0.0512$ & $30.9750$\\
$\text{Optimal-3}$ & $8.0043$ & $0.0143$ & $58.5620$\\
\bottomrule
\end{tabular}
\end{table}

Fig.\ref{fig:2} shows that the PoEMV-1, PoEMV-2 and CoEMV algorithms effectively bring the terminal net wealth to approach the target value $d=8$ within a finite number of iterations, which is consistent with the policy convergence results proved in Section \ref{subsec 4.1},
while their Lagrange multipliers exhibit similar patterns and stabilize near \(8\).
Table \ref{Table 1} demonstrates that the three learned policies have out-of-sample means between \(7.9954\) and \(8.0074\), with comparable variances. PoEMV-1 has the lowest variance and the highest Sharpe ratio among the learning-based methods. Among the model-based references, all three attain means closer to \(d=8\), and the posterior-projection policy (Optimal-1) reduces the variance of terminal net wealth and markedly improves the Sharpe ratio relative to the no-learning reference (Optimal-2), highlighting the value of dynamically updating market-state information in multi-period portfolio decisions.
Across the three pairwise comparisons in Table \ref{Table 1}, in which the realized means differ only slightly, PoEMV-1 and CoEMV exhibit modest increases in out-of-sample variance and declines in the Sharpe ratio relative to their respective model-based references, whereas PoEMV-2 achieves a lower variance and a higher Sharpe ratio. By maintaining low terminal net-wealth variance while attaining a mean close to \(d=8\), PoEMV-1 effectively captures the optimal policy and thus offers a viable alternative to model-based optimization.

\subsection{Empirical study}

In this section, we assess the performance of PoEMV-1 using real financial data. The initial market state distribution $\hat{p}_0$ is set to \(\hat{p}_0=0.3\), as in Section~\ref{simu}. Conditional on the historical market regime, the one-step liability returns \(q_t(1)-1\) and \(q_t(2)-1\) are simulated from \(\mathcal{N}(0.05\Delta t,0.10\Delta t)\) and \(\mathcal{N}(0.01\Delta t,0.20\Delta t)\), respectively, where \(\Delta t=1/252\) for daily rebalancing and \(\Delta t=1/12\) for monthly rebalancing. The return rate of the risk-free asset \(e_t^0\) is represented by a common observable cash rate estimated from the pooled training-period daily 13-week U.S. Treasury-bill yields (\texttt{\textasciicircum{}IRX}) and is kept identical across market regimes and competing methods. For the risky asset, we use the S\&P 500 index and collect 30 years of daily data from March 1, 1995, to May 30, 2025, sourced from Yahoo Finance and CRSP. Monthly returns are obtained by compounding the daily returns within each calendar month.
Our dataset includes 289 stocks with complete 30-year records. The first 20 years of data are used for training and validation, while the remaining 10 years are allocated for out-of-sample testing.

\begin{figure}[htbp]
	\centering
	\includegraphics[width=5in]{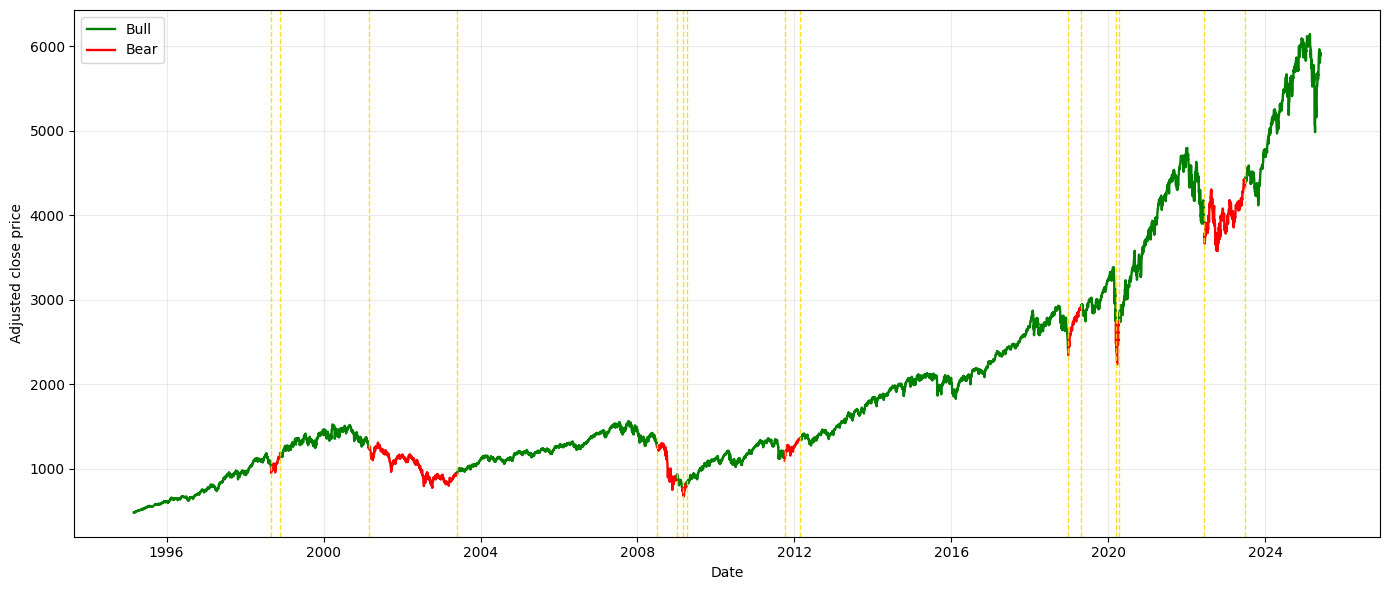}
	\caption{S\&P 500 index from 01/03/1995 to 30/05/2025 with bull (green) and bear (red) market phases.}
	\label{fig:5}
\end{figure}

Figure \ref{fig:5} illustrates the regime-switching behavior of the S\&P 500 index during this period, revealing the alternation between bull (green) and bear (red) market phases. A bull market is defined by an increase of at least $\gamma_1=24\%$, while a bear market is characterized by a decline of at least $\gamma_2=19\%$. This regime classification provides a clear view of the market's cyclical nature and sets the stage for evaluating PoEMV-1's effectiveness in navigating these transitions.

Our objective is to construct a portfolio that dynamically allocates between risk-free and risky assets for the mean-variance (MV) problem over a 10-year investment horizon. The terminal target is set at $d=8$, with initial wealth $x_0=1$ and initial liability $l_0=0.1$, corresponding to a 24.42\% annualized growth rate in net wealth. Unlike in the simulation study, where unlimited paths can be generated, we face the challenge of data insufficiency with only one observed trajectory for the index. To address this, we use a panel of 79 stocks whose market betas, estimated from the first 20 years of data, are between 0.9 and 1.1. We then apply the universal training method proposed by \cite{Wang2020}, treating the 79 paths of different stocks as those of a single stock. Within the training period, 63 stocks are assigned to policy fitting and 16 are held out for validation. The same split is used at both frequencies.

For each of the 63 fitting stocks, the first 20 years of data are divided into overlapping 10-year blocks. The \(5{,}034\) daily and \(239\) monthly training observations per stock yield \(158{,}445\) $(=(5{,}034 - 10 \times 252 + 1) \times 63)$ daily and \(7{,}560\) $(=((20\times12-1)-10\times12+1)\times63)$) monthly fitting blocks, respectively.
At each frequency, the market parameters are estimated from the training sample. Specifically, for regime \(i\), let
\(r_{i,1},\ldots,r_{i,n_i}\) denote the one-period net returns pooled across the 63 fitting stocks for intervals that begin in regime \(i\).
The mean \(\mathbb{E}[e_t^1(i)]\) and variance \(\mathrm{Var}[e_t^1(i)]\) of the gross risky-asset return are estimated by \(1+\bar r_i\) and \(s_i^2\), respectively, where
$\bar r_i=\sum_{k=1}^{n_i}r_{i,k}/n_i$ and
$s_i^2=\sum_{k=1}^{n_i}(r_{i,k}-\bar r_i)^2/(n_i-1)$.
Transition probabilities are estimated as $\widehat P_{ij}=N_{ij}/(N_{i1}+N_{i2})$, $i,j\in\{1,2\}$, where \(N_{ij}\) counts transitions from regime \(i\) to regime \(j\) between consecutive labels in the training-period market regime sequence. These estimates are held fixed throughout training.
Each sampled 10-year block supplies the observed stock returns for policy evaluation and policy-gradient updates. At the end of each period, PoEMV-1 updates its posterior belief using the observed stock return and the simulated liability return. The learning rates $\alpha$ and $\eta$ match those used in the simulation study.

We compare two algorithms: our PoEMV-1 algorithm and the EMV algorithm from \cite{Cui2023}. Both algorithms parameterize the value functions and policies using semi-analytical representations, followed by actor-critic learning with parameters updated via the martingale loss (ML) and the policy gradient. Each algorithm is run for \(20{,}000\) iterations at each frequency, drawing eight blocks per iteration. The policies are selected using \(256\) fixed validation blocks from the held-out stocks, with a mean-target tolerance of \(0.05\) and the same selection rule as in the simulation study.

The EMV policy in \cite{Cui2023} is derived from a multi-period model with one risk-free asset and one risky asset, without regime-switching, so no filtering is required. In the empirical comparison, EMV uses the same stock returns, cash rate and liability specification as PoEMV-1, with the belief coordinate set to zero.
To facilitate comparative analysis, we also implement a model-based belief-grid reference. Its policy is obtained by backward recursion on a grid of posterior beliefs using the fitted market parameters, and its multiplier is calibrated on \(160\) fitting blocks to match the target mean.

For out-of-sample testing, we apply the selected PoEMV-1 and EMV policies and the belief-grid reference to the same 79 stock paths at each frequency. The daily evaluation uses the first \(2{,}520\) testing returns (March 2, 2015, to March 6, 2025), and the monthly evaluation uses the first \(120\) testing returns (March 2015 to February 2025). As in the simulation study, actions are sampled from the Gaussian policies and clipped to \([-20,20]\), with the Gaussian variance constrained to \([10^{-8},5\times10^{-4}]\) in the empirical study. The liability and policy shocks are paired across methods. Table \ref{Table 2} presents the mean and variance of terminal net wealth, along with the 10-year Sharpe ratio, for these implementations.

\begin{table}[!ht]
\caption{Out-of-sample performance under daily and monthly rebalancing.}
\label{Table 2}
\centering
\begin{tabular}{l c c c}
\toprule
$\ $ & $\text{Mean}$ & $\text{Variance}$ & $\text{Sharpe ratio}$ \\
\cline{1-4}
$\text{PoEMV-1}$ & $7.9194$ & $0.0337$ & $37.6756$\\
$\text{EMV}$ & $7.9173$ & $0.0520$ & $30.3451$ \\
$\text{Classical}$ & $7.9576$ & $0.9330$ & $7.2030$\\
\bottomrule
\multicolumn{2}{c}{Daily rebalancing\ \ \ \ \ }\\
\toprule
$\ $ & $\text{Mean}$ & $\text{Variance}$ & $\text{Sharpe ratio}$ \\
\cline{1-4}
$\text{PoEMV-1}$ & $7.9992$ & $0.0811$ & $24.5745$\\
$\text{EMV}$ & $7.9857$ & $0.0818$ & $24.4243$ \\
$\text{Classical}$ & $8.4777$ & $1.1980$ & $6.8319$\\
\bottomrule
\multicolumn{2}{c}{Monthly rebalancing}
\end{tabular}
\end{table}

As shown in Table \ref{Table 2}, under daily rebalancing all three policies produce mean terminal net wealth close to $d=8$. Under monthly rebalancing, PoEMV-1 and EMV remain close to the target, while the classical method attains a mean of $8.4777$. A significant differences emerge in the variances and Sharpe ratios across methods. At both rebalancing frequencies, PoEMV-1 achieves the largest Sharpe ratio, more than three times that of the belief-grid reference, which exhibits the largest variance in terminal net wealth, largely due to extremely low outcomes observed in certain stock paths.
The comparison between PoEMV-1 and EMV underscores the critical importance of incorporating market regime filtering as a guiding signal for long-term investment decisions. The Sharpe ratio of the EMV algorithm is lower than that of the completely observable COEMV (see Table \ref{Table 1}), highlighting the practical gains made by accounting for market regimes in optimal portfolio construction. Moreover, while both the classical method and PoEMV-1 utilize the same estimated environment parameters for filtering, their approach to deriving the optimal policy differs fundamentally. The classical method relies heavily on the estimated parameters, which can introduce bias and inaccuracies into the optimal policy. In contrast, PoEMV-1 mitigates this risk by refining the policy through iterative interactions, leading to superior outcomes, particularly in environments where estimation errors are prevalent.

Additionally, the impact of rebalancing frequency can be seen when comparing daily and monthly results. Across all three algorithms, the variance of terminal wealth is generally lower under daily rebalancing than under monthly rebalancing. This suggests that implementing the learned policy on a daily basis enables more timely adjustments to the investment strategy in response to market changes, resulting in better overall investment performance.

\section{Conclusion}\label{sec6}

In this study, we modeled the optimal portfolio problem under the exploratory mean-variance (EMV) criterion by incorporating assets, liabilities, and market regimes, with the latter governed by a multi-period homogeneous Markov chain to reflect realistic financial market dynamics. Our key findings highlight that, under partial information-where the investor can only observe stock prices and not the underlying market regime-the original problem can be effectively transformed into a fully observable investment optimization problem using stochastic filtering theory. This transformation enabled us to derive the value function and the optimal investment strategy based on the Bellman principle.

A significant contribution of this research is the development of actor-critic reinforcement learning (RL) algorithms tailored to multi-period EMV investment in a regime-switching market. Specifically, the PoEMV-1, PoEMV-2, and CoEMV algorithms were introduced to capture optimal investment policies under different information settings. Our findings demonstrate that PoEMV-1 consistently achieves superior empirical performance and captures the optimal policy with high accuracy compared to CoEMV and PoEMV-2. Notably, PoEMV-1 outperforms the classical stochastic control approach and the standard EMV algorithm, which do not account for market regime changes. This result emphasizes the importance of incorporating regime-switching dynamics in portfolio optimization to enhance decision-making and performance.

Our analysis provides a robust economic interpretation of these comparative results, further validating the practical advantages of using regime-switching models in real-world financial settings. Looking forward, future research should focus on advancing RL algorithms for investment and trading problems, particularly by incorporating regime-switching dynamics alongside alternative performance metrics, to further improve decision-making frameworks in complex financial environments.
These findings underscore the potential of regime-aware approaches in portfolio optimization, highlighting the significant benefits they offer over more traditional models.

\section*{Conflict of interest}
All authors declare no conflicts of interest in this paper.







\appendix

\section{Proofs and more}\label{Section appen}

In the following we provide a complete proof of the theorems presented in this paper.
\subsection{Proof of Theorem \ref{thm32a}}\label{appen11}
\begin{proof}
The transformation from $(e_t^0,e_t^1)$ to $(e_t^0,e_t^1-e_t^0)$ is invertible. Positive definiteness thus gives $B_{1,t}^{(i)}>0$ and $F_{1,t}^{(i)}>0$. The recursion from $\bar F_{1,T}^{(i)}=1$ yields $\bar F_{1,t}^{(i)}>0$ for all $t$ and $i$, so each action curvature in \eqref{eq:general-b01} is positive. 
For the regime-switching CoEMV problem \eqref{max66}, the current regime $\varepsilon_t=i$ is known at time $t$. For each fixed $w\in\mathbb{R}$, by the Bellman's principle, we have
\begin{align}\label{jtxp1}
&{J}_t^{c*}(x,l,i; w)
= \min_{\pi_t^c} \int_{\mathbb{R}} \pi_t^c(u) \, \Big(\sum_{j\in \mathcal{M}} P_{ij}\, \mathbb{E}_{t,x,l,i}\big[{J}_{t+1}^{c*}({x}_{t+1}^{\pi^c},{l}_{t+1},j; w) \mid u_t^{\pi^c}=u\big] + \lambda \ln\pi_t^{c}(u)\Big) \mathrm{d} u.
\end{align}
We now consider the problem \eqref{max66} from period $T$ to $T-1$.
Since $J_{T}^{c*}$ in \eqref{term098} is independent of the next-period regime $\varepsilon_T=j\in \mathcal{M}$ and $\sum_{j\in\mathcal{M}} P_{ij}=1$, we have $\sum_{j\in\mathcal{M}} P_{ij}J_T^{c*}(x, l, j;w)=(x-l-w)^2 - (w-d)^2$. By \eqref{xt2}, \eqref{lt} and \eqref{term098}, we get
\begin{align*}
&{J}_{T-1}^{c*}(x,l,i; w)
=\min_{\pi_{T-1}^c} \int_{\mathbb{R}}\pi^c_{T-1}(u)
\Big( \mathbb{E}_{T-1,x,l,i}
\big[({x}_T^{\pi^c}-{l}_T-w)^2 \mid u_{T-1}^{\pi^c}=u\big]
+\lambda\ln\pi^c_{T-1}(u) \Big) \mathrm{d} u-(w-d)^2\\
&=\min_{\pi_{T-1}^c}  \mathbb{E}_{T-1,x,l,i}
\big[({x}_T^{\pi^c})^2
-2{x}_T^{\pi^c}{l}_T-2w{x}_T^{\pi^c}+2w{l}_T+{l}_T^2
+\lambda\int_{\mathbb{R}}\pi^c_{T-1}(u)\ln\pi^c_{T-1}(u)\mathrm{d}u \big]-d^2+2wd\\
&= \min_{\pi_{T-1}^c}  \mathbb{E}_{T-1,x,l,i}
\big[\big(\zeta_{0, T-1}^{(i)} {x}_{T-1}^{\pi^c}\big)^2
+\big(\zeta_{1, T-1}^{(i)}
u_{T-1}^{\pi^c}\big)^2
+2\zeta_{0, T-1}^{(i)} {x}_{T-1}^{\pi^c}
\zeta_{1, T-1}^{(i)}u_{T-1}^{\pi^c}
\\
& \qquad \qquad \qquad \qquad
-2\zeta_{0, T-1}^{(i)} {x}_{T-1}^{\pi^c}\zeta_{2, T-1}^{(i)} {l}_{T-1}
-2\zeta_{1, T-1}^{(i)} u_{T-1}^{\pi^c} \zeta_{2, T-1}^{(i)} {l}_{T-1}
-2w\zeta_{0, T-1}^{(i)} {x}_{T-1}^{\pi^c}
-2w\zeta_{1, T-1}^{(i)}u_{T-1}^{\pi^c}
\\
&\qquad \qquad \qquad \qquad
+2w\zeta_{2, T-1}^{(i)} {l}_{T-1}
+\big(\zeta_{2, T-1}^{(i)} {l}_{T-1}\big)^2
+\lambda\int_{\mathbb{R}}\pi^c_{T-1}(u)\ln\pi^c_{T-1}(u)\mathrm{d}u
\big]-d^2+2wd\\
&= \min_{\pi_{T-1}^c} \int_{\mathbb{R}}
\Big(
\mathbb{E}[(\zeta_{1, T-1}^{(i)})^{2}] u^2
+2\mathbb{E}[\zeta_{0, T-1}^{(i)} \zeta_{1, T-1}^{(i)}] x u
-2w\mathbb{E}[\zeta_{1, T-1}^{(i)}] u
-2\mathbb{E}[\zeta_{1, T-1}^{(i)} \zeta_{2, T-1}^{(i)}]l u
+\lambda \ln\pi^c_{T-1}(u)
\Big)\pi^c_{T-1}(u) \mathrm{d}u \\
& \qquad \qquad
+\mathbb{E}[(\zeta_{0, T-1}^{(i)}){^2}]x^2-2\mathbb{E}[\zeta_{0, T-1}^{(i)} \zeta_{2, T-1}^{(i)}] l x
-2w\mathbb{E}[\zeta_{0, T-1}^{(i)}] x
+2w\mathbb{E}[\zeta_{2, T-1}^{(i)}] l
+\mathbb{E}[(\zeta_{2, T-1}^{(i)})^2] l^2
-d^2+2wd \\
&=\min_{\pi_{T-1}^c}  \int_{\mathbb{R}}
\Big(
{B}_{1,T-1}^{(i)} u^2
+2\mu_{T-1}^{(i)}(x,l;w) u
+\lambda \ln\pi^c_{T-1}(u)
\Big)\pi^c_{T-1}(u) \mathrm{d}u\\
&\qquad \qquad +{B}_{0,T-1}^{(i)}x^2
-2{C}_{02,T-1}^{(i)} l x-2w{A}_{0,T-1}^{(i)} x
+2w{A}_{2,T-1}^{(i)} l
+{B}_{2, T-1}^{(i)} l^2
-d^2+2wd,
\end{align*}
where $\mu_{T-1}^{(i)}(x,l;w):=C_{01, T-1}^{(i)} x -w {A}_{1,T-1}^{(i)} - C_{12,T-1}^{(i)} l$.
In order to minimize the above function with respect to $\pi^c_{T-1}$, we consider the variational first-order condition. This yields
\begin{align*}
{B}_{1,T-1}^{(i)} u^2
+2\mu_{T-1}^{(i)}(x,l;w)  u
+\lambda \ln\pi^c_{T-1}(u)+\lambda+\kappa=0.
\end{align*}
The multiplier $\kappa$ enforces $\int_{\mathbb R}\pi^c_{T-1}(u)\,\mathrm du=1$. Since $B_{1,T-1}^{(i)}>0$, the resulting normalized density is strictly positive and minimizes the stage objective by the Gibbs variational identity. Its density is
\begin{align*}
\pi_{T-1}^{c*}(u; x, l, i,w)
&=\frac{\exp\big\{-\frac{1}{\lambda}\big({B}_{1,T-1}^{(i)} u^2
+2\mu_{T-1}^{(i)}(x,l;w)  u\big)\big\}}
{\int_{\mathbb{R}}\exp\{-\frac{1}{\lambda}\big({B}_{1,T-1}^{(i)} u^2
+2\mu_{T-1}^{(i)}(x,l;w) u\big)\}\mathrm{d} u}
=\mathcal{N}\Big(u\Big{|}-\frac{\mu_{T-1}^{(i)}(x,l;w)}{{B}_{1,T-1}^{(i)}}, \frac{\lambda}{2 {B}_{1,T-1}^{(i)}}\Big).
\end{align*}
Then, we have
\begin{align*}
&\int_{\mathbb{R}} \ln \pi_{T-1}^{c*}(u) \pi_{T-1}^{c*}(u) \mathrm{d} u
=\int_{\mathbb{R}}\ln\Big(\sqrt{\frac{{B}_{1,T-1}^{(i)}}{\pi\lambda}}
\exp \Big\{-\frac{{B}_{1,T-1}^{(i)}}{\lambda}\big(u +\frac{\mu_{T-1}^{(i)}(x,l;w)}{{B}_{1,T-1}^{(i)}}\big)^2 \Big\}\Big) \pi_{T-1}^{c*}(u) \mathrm{d} u\\
&=\frac{1}{2}\ln \frac{{B}_{1,T-1}^{(i)}}{\pi\lambda}-\frac{{B}_{1,T-1}^{(i)}}{\lambda}
\big((\frac{\mu_{T-1}^{(i)}(x,l;w)}{{B}_{1,T-1}^{(i)}})^2+\frac{\lambda}{2 {B}_{1,T-1}^{(i)}} -2(\frac{\mu_{T-1}^{(i)}(x,l;w)}{{B}_{1,T-1}^{(i)}})^2 +(\frac{\mu_{T-1}^{(i)}(x,l;w)}{{B}_{1,T-1}^{(i)}})^2\big)\\
&=\frac{1}{2}\ln \frac{{B}_{1,T-1}^{(i)}}{\pi\lambda} -\frac{1}{2}.
\end{align*}
Substituting $\pi_{T-1}^{c*}(u; x, l, i,w)$ and the above result into ${J}_{T-1}^{\pi^c}(x, l, i; w)$ leads to
\begin{align*}
&{J}_{T-1}^{c*}(x, l, i; w)
=
{B}_{1,T-1}^{(i)} \Big(\big(\frac{\mu_{T-1}^{(i)}(x,l;w)}{{B}_{1,T-1}^{(i)}}\big)^2+\frac{\lambda}{2 {B}_{1,T-1}^{(i)}}\Big)
-\frac{2\big(\mu_{T-1}^{(i)}(x,l;w)\big)^2}{{B}_{1,T-1}^{(i)}}
+\frac{\lambda}{2}\ln \frac{{B}_{1,T-1}^{(i)}}{\pi\lambda} -\frac{\lambda}{2}
+{B}_{0,T-1}^{(i)}x^2\\
&\qquad \qquad \qquad \quad \
-2{C}_{02,T-1}^{(i)} l x
-2w{A}_{0,T-1}^{(i)} x
+2w{A}_{2,T-1}^{(i)} l
+{B}_{2,T-1}^{(i)} l^2
 -d^2+2wd\\
&=\frac{\lambda}{2}\ln \frac{{B}_{1,T-1}^{(i)}}{\pi\lambda}
-\frac{\big(\mu_{T-1}^{(i)}(x,l;w)\big)^2}{{B}_{1,T-1}^{(i)}}
+{B}_{0,T-1}^{(i)}x^2
-2(w{A}_{0,T-1}^{(i)}+{C}_{02,T-1}^{(i)} l) x
+2w{A}_{2,T-1}^{(i)} l
 +{B}_{2,T-1}^{(i)} l^2
-d^2+2wd\\
&=\frac{\lambda}{2}\ln \frac{{B}_{1,T-1}^{(i)}}{\pi\lambda}
+\frac{{F}_{1,T-1}^{(i)}}{{B}_{1,T-1}^{(i)}} x^2
-\frac{2 (w {F}_{2,T-1}^{(i)} +
 {F}_{3,T-1}^{(i)}l )}{{B}_{1,T-1}^{(i)}}x
-\frac{(w{A}_{1,T-1}^{(i)}+{C}_{12,T-1}^{(i)}l)^2}{{B}_{1,T-1}^{(i)}}
+2{A}_{2,T-1}^{(i)} wl +{B}_{2,T-1}^{(i)} l^2
-d^2+2wd,
\end{align*}
where $F_{r, t}^{(i)}$ for $r=1,2,3$ are defined by \eqref{f123t}.

We then consider the EMV problem \eqref{max66} from period $T-1$ to $T-2$. Note that the value function ${J}_{T-1}^{c*}$ depends on the next-period regime $\varepsilon_{T-1}=j\in \mathcal M$, unlike the terminal step. This implies
\begin{align}\label{eq:A1-coupled-Tm1}
\sum_{j\in\mathcal M}P_{ij}J_{T-1}^{c*}(x',l',j;w)
&= \frac{\lambda}{2}\hat L_{T-1}^{(i)} +
\hat F_{1,T-1}^{(i)}(x')^2
-2\bigl(w\hat F_{2,T-1}^{(i)}+\hat F_{3,T-1}^{(i)} l'\bigr)x'
-w^2\hat D_{1,T-1}^{(i)}
+2wl'\bigl(\hat A_{2,T-1}^{(i)} -\hat D_{2,T-1}^{(i)}\bigr)\nonumber \\
& \quad + (l')^2\bigl(\hat B_{2,T-1}^{(i)} - \hat D_{3,T-1}^{(i)}\bigr)
-d^2+2wd,
\end{align}
where the transition-weighted coefficients are defined by, for $r=1,2, 3$,
\begin{align*}
&
\hat F_{r,T-1}^{(i)}
:= \sum_{j\in\mathcal M}P_{ij}\frac{F_{r,T-1}^{(j)}}{B_{1,T-1}^{(j)}}, \quad
\hat L_{T-1}^{(i)}
:=
\sum_{j\in\mathcal M}P_{ij}\ln\frac{B_{1,T-1}^{(j)}}{\pi\lambda},
\quad
\hat A_{2,T-1}^{(i)}
:= \sum_{j\in\mathcal M}P_{ij}A_{2,T-1}^{(j)},
\quad
\hat B_{2,T-1}^{(i)}
:= \sum_{j\in\mathcal M}P_{ij}B_{2,T-1}^{(j)}, \\
&
\hat D_{1,T-1}^{(i)}
:=
\sum_{j\in\mathcal M}P_{ij}\frac{\bigl(A_{1,T-1}^{(j)}\bigr)^2}{B_{1,T-1}^{(j)}},\quad
\hat D_{2,T-1}^{(i)}
:=
\sum_{j\in\mathcal M}P_{ij}\frac{A_{1,T-1}^{(j)}C_{12,T-1}^{(j)}}{B_{1,T-1}^{(j)}},\quad
\hat D_{3,T-1}^{(i)}
:=
\sum_{j\in\mathcal M}P_{ij}\frac{\bigl(C_{12,T-1}^{(j)}\bigr)^2}{B_{1,T-1}^{(j)}}.
\end{align*}
Then, at period $T-2$, the Bellman recursion \eqref{jtxp1} becomes
\begin{align*}
&{J}_{T-2}^{c*}(x, l, i; w)
= \min_{\pi_{T-2}^c} \int_{\mathbb{R}} \pi_{T-2}^c(u) \, \Big( \mathbb{E}_{T-2, x, l, i} \big[\sum_{j\in \mathcal{M}} P_{ij}\, {J}_{T-1}^{c*}({x}_{T-1}^{\pi^c},{l}_{T-1},j; w) \mid u_{T-2}^{\pi^c}=u\big] + \lambda \ln\pi_{T-2}^{c}(u)\Big) \mathrm{d} u \\
&=\min_{\pi_{T-2}^c} \int_{\mathbb{R}} \pi_{T-2}^c(u) \,  \Big( \mathbb{E}_{T-2,x, l, i}
\big[\hat F_{1,T-1}^{(i)}(x_{T-1}^{\pi^c})^2
-2\bigl(w\hat F_{2,T-1}^{(i)}+\hat F_{3,T-1}^{(i)} l_{T-1}\bigr)x_{T-1}^{\pi^c}
-w^2\hat D_{1,T-1}^{(i)} \nonumber \\
& \quad
+ 2wl_{T-1} \bigl(\hat A_{2,T-1}^{(i)}-\hat D_{2,T-1}^{(i)}\bigr) + (l_{T-1})^2\bigl(\hat B_{2,T-1}^{(i)} -\hat D_{3,T-1}^{(i)}\bigr)\big]
+\lambda \ln\pi^c_{T-2}(u) \Big) \mathrm{d} u
+\frac{\lambda}{2}\hat L_{T-1}^{(i)} -d^2+2wd.
\end{align*}
Substituting the wealth process \eqref{xt2} and the liability process \eqref{lt} into this equation yields
\begin{align*}
&{J}_{T-2}^{c*}(x, l, i; w)
= \min_{\pi_{T-2}^c} \int_{\mathbb{R}} \pi_{T-2}^c(u) \, \Big(
\mathbb{E}_{T-2,x, l, i}
\big[\hat{F}_{1,T-1}^{(i)}
(\zeta_{0, T-2}^{(i)} x+\zeta_{1, T-2}^{(i)} u)^2\\
&\ \ \
-2(w{\hat F}_{2,T-1}^{(i)} +{\hat F}_{3,T-1}^{(i)} {\zeta}_{2, T-2}^{(i)} {l})
(\zeta_{0, T-2}^{(i)} {x} + \zeta_{1, T-2}^{(i)} u)
- w^2\hat D_{1,T-1}^{(i)}
+ 2w \zeta_{2, T-2}^{(i)} l \bigl(\hat A_{2,T-1}^{(i)} -\hat D_{2,T-1}^{(i)}\bigr) \nonumber \\
& \quad
+(\zeta_{2, T-2}^{(i)} l)^2\bigl(\hat B_{2,T-1}^{(i)} - \hat D_{3,T-1}^{(i)}\bigr) \big]
+\lambda \ln\pi^c_{T-2}(u) \Big) \mathrm{d} u + \frac{\lambda}{2}\hat L_{T-1}^{(i)}
-d^2+2wd \\
&=
\min_{\pi_{T-2}^c}
\int_{\mathbb R}\pi_{T-2}^c(u)
\Bigl[
\bar B_{1,T-2}^{(i)}u^2+2\mu_{T-2}^{(i)}(x,l;w)u+\lambda\ln\pi_{T-2}^c(u)
\Bigr] \mathrm du
+\Xi_{T-2}^{(i)}(x,l;w),
\end{align*}
where
\begin{align*}
& \bar B_{1,T-2}^{(i)}
:= \hat F_{1,T-1}^{(i)}\,B_{1,T-2}^{(i)},\\
& \mu_{T-2}^{(i)}(x,l;w)
:= \hat F_{1,T-1}^{(i)}C_{01,T-2}^{(i)}x
-w\hat F_{2,T-1}^{(i)}A_{1,T-2}^{(i)}
-l\hat F_{3,T-1}^{(i)}C_{12,T-2}^{(i)},\\
&\Xi_{T-2}^{(i)}(x,l;w)
:= \hat F_{1,T-1}^{(i)}B_{0,T-2}^{(i)}x^2
-2\Bigl( w\hat F_{2,T-1}^{(i)}A_{0,T-2}^{(i)}
+
l\hat F_{3,T-1}^{(i)}C_{02,T-2}^{(i)}
\Bigr)x -w^2\hat D_{1,T-1}^{(i)} \\
&\qquad \qquad \qquad \quad
+ 2wl\,A_{2,T-2}^{(i)}\bigl(\hat A_{2,T-1}^{(i)} - \hat D_{2,T-1}^{(i)}\bigr)
+ l^2B_{2,T-2}^{(i)}\bigl(\hat B_{2,T-1}^{(i)} - \hat D_{3,T-1}^{(i)}\bigr)
+\frac{\lambda}{2}\hat L_{T-1}^{(i)}
-d^2+2wd.
\end{align*}
By the same variational argument as before, we obtain an optimal feedback control whose density function is
\begin{align*}
&\pi_{T-2}^{c*}(u; x, l, i,w)
=\mathcal{N}\Big(u\Big{|}-\frac{\mu_{T-2}^{(i)}(x,l;w)}{{\bar B}_{1,T-2}^{(i)}},
\frac{\lambda}{2 {\bar B}_{1,T-2}^{(i)}}\Big).
\end{align*}
Then, we have $\int_{\mathbb{R}} \ln \pi_{T-2}^{c*}(u) \pi_{T-2}^{c*}(u) \mathrm{d} u
=\frac{1}{2}\ln\frac{{\bar B}_{1,T-2}^{(i)}}{\pi \lambda}
-\frac{1}{2}$.
Substituting $\pi_{T-2}^{c*}(u; x, l, i,w)$ and the above result into ${J}_{T-2}^{c*}(x, l, i; w)$, we arrive at
\begin{align*}
&{J}_{T-2}^{c*}(x, l, i; w)
=
-\frac{\bigl(\mu_{T-2}^{(i)}(x,l;w)\bigr)^2}{\bar B_{1,T-2}^{(i)}}
+\frac{\lambda}{2}\ln\frac{\bar B_{1,T-2}^{(i)}}{\pi\lambda}
+\Xi_{T-2}^{(i)}(x,l;w).
\end{align*}
Expanding the square term and simplifying, we obtain
\begin{align}\label{eq:J-T-2-compact}
J_{T-2}^{c*}(x,l,i;w)
&=
\frac{\lambda}{2}\bar L_{T-2}^{(i)}
+\bar F_{1,T-2}^{(i)}x^2
-2\bigl(w\bar F_{2,T-2}^{(i)}+l\bar F_{3,T-2}^{(i)}\bigr)x
-w^2\bar D_{1,T-2}^{(i)}
+2wl\bigl(\bar A_{2,T-2}^{(i)} -\bar D_{2,T-2}^{(i)}\bigr) \nonumber \\
& \quad + l^2\bigl(\bar B_{2,T-2}^{(i)} - \bar D_{3,T-2}^{(i)}\bigr)
-d^2+2wd,
\end{align}
where, for $r=1,2, 3$,
\begin{align*}
& \bar F_{r,T-2}^{(i)}
:=
\hat F_{r,T-1}^{(i)}\frac{F_{r,T-2}^{(i)}}{B_{1,T-2}^{(i)}}, \quad
\bar L_{T-2}^{(i)} := \ln\frac{\bar B_{1,T-2}^{(i)}}{\pi\lambda}
+\hat L_{T-1}^{(i)},
\quad
\bar A_{2,T-2}^{(i)}
:=
A_{2,T-2}^{(i)}\hat A_{2,T-1}^{(i)},
\quad
\bar B_{2,T-2}^{(i)} :=
B_{2,T-2}^{(i)}\hat B_{2,T-1}^{(i)}, \\
& \bar D_{1,T-2}^{(i)}
:=
\hat D_{1,T-1}^{(i)}
+\frac{\bigl(\hat F_{2,T-1}^{(i)}A_{1,T-2}^{(i)}\bigr)^2}
{\bar B_{1,T-2}^{(i)}}, \qquad
\bar D_{2,T-2}^{(i)}
:=
A_{2,T-2}^{(i)}\hat D_{2,T-1}^{(i)}
+\frac{\hat F_{2,T-1}^{(i)}\hat F_{3,T-1}^{(i)}A_{1,T-2}^{(i)}C_{12,T-2}^{(i)}}
{\bar B_{1,T-2}^{(i)}}, \\
& \bar D_{3,T-2}^{(i)}
:=
B_{2,T-2}^{(i)}\hat D_{3,T-1}^{(i)}
+\frac{\bigl(\hat F_{3,T-1}^{(i)}C_{12,T-2}^{(i)}\bigr)^2}
{\bar B_{1,T-2}^{(i)}}.
\end{align*}

Next, we state the induction hypothesis in a form that is consistent with the
backward derivations from $T$ to $T-1$, and from $T-1$ to $T-2$. For each fixed $w\in\mathbb R$, we claim that, for every $t\in\{0,1,\ldots,T-1\}$ and
$j\in\mathcal M$, the value function at period $t+1$ admits
\begin{align}\label{eq:IH-value-t-plus-one}
J_{t+1}^{c*}(x,l,j;w)
&=
\frac{\lambda}{2}\bar L_{t+1}^{(j)}
+\bar F_{1,t+1}^{(j)}x^2
-2\bigl(w\bar F_{2,t+1}^{(j)}
+l\bar F_{3,t+1}^{(j)}\bigr)x
\notag\\
&\quad
-\bar D_{1,t+1}^{(j)}w^2
+2\bigl(\bar A_{2,t+1}^{(j)}
-\bar D_{2,t+1}^{(j)}\bigr)wl
+\bigl(\bar B_{2,t+1}^{(j)}
-\bar D_{3,t+1}^{(j)}\bigr)l^2
-d^2+2wd,
\end{align}
where the associated coefficients are given by \eqref{f123t}-\eqref{eq:barD3-general}.
Then, by the Bellman principle, under the current regime $\varepsilon_t=i$, we have
\begin{align*}
J_t^{c*}(x,l,i;w)
&=
\min_{\pi_t^c}
\int_{\mathbb R}\pi_t^c(u)
\Big(
\mathbb E_{t, x, l, i}\big[
\sum_{j\in\mathcal M}P_{ij} J_{t+1}^{c*}(x_{t+1}^{\pi^c},l_{t+1},j;w) \mid u_t^{\pi^c}=u \big]
+\lambda\ln\pi_t^c(u)
\Big)\mathrm du\\
&=
\min_{\pi_t^c}
\int_{\mathbb R}\pi_t^c(u)
\Big[
\bar B_{1,t}^{(i)}u^2
+2\mu_t^{(i)}(x,l;w)u
+\lambda\ln\pi_t^c(u)
\Big] \mathrm du
+\Xi_t^{(i)}(x,l;w),
\end{align*}
where $\bar B_{1,t}^{(i)}$ and $\mu_t^{(i)}(x,l;w)$ are defined by \eqref{eq:general-b01} and \eqref{eq:general-mu}, respectively, and
\begin{align*}
&\Xi_t^{(i)}(x,l;w) =
B_{0,t}^{(i)}x^2
\sum_{j\in\mathcal M}P_{ij}\bar F_{1,t+1}^{(j)}
-2\Big(
wA_{0,t}^{(i)}
\sum_{j\in\mathcal M}P_{ij}\bar F_{2,t+1}^{(j)}
+
lC_{02,t}^{(i)}
\sum_{j\in\mathcal M}P_{ij}\bar F_{3,t+1}^{(j)}
\Big)x
-w^2 \sum_{j\in\mathcal M}P_{ij}\bar D_{1,t+1}^{(j)} \\
&
+2wl A_{2,t}^{(i)}
\Big(
\sum_{j\in\mathcal M}P_{ij}\bar A_{2,t+1}^{(j)}
-
\sum_{j\in\mathcal M}P_{ij}\bar D_{2,t+1}^{(j)}
\Big)
+l^2B_{2,t}^{(i)}
\Big(
\sum_{j\in\mathcal M}P_{ij}\bar B_{2,t+1}^{(j)}
-
\sum_{j\in\mathcal M}P_{ij}\bar D_{3,t+1}^{(j)}
\Big)
+\frac{\lambda}{2}
\sum_{j\in\mathcal M}P_{ij}\bar L_{t+1}^{(j)}
-d^2+2wd.
\end{align*}
By the same variational argument as before, we have $\pi_t^{c*}(u;x,l,i,w)$ satisfies \eqref{pi32a}. Equivalently, $\int_{\mathbb R}\ln\pi_t^{c*}(u)\pi_t^{c*}(u)\,\mathrm du
=
\frac12\ln\frac{\bar B_{1,t}^{(i)}}{\pi\lambda}
-\frac12$, and then $J_t^{c*}(x,l,i;w)$ has the form of \eqref{eq:general-value-final0}. Therefore, by the principle of mathematical induction, for all $t\in\{0,1,\ldots,T\}$ and all $i\in\mathcal M$, the value function of the CoEMV problem \eqref{max66} at period $t$ has the form of \eqref{eq:general-value-final0}, while the optimal relaxed feedback control is Gaussian with density \eqref{pi32a}.

\end{proof}

\subsection{Proof of Theorem \ref{them1122}}\label{appen1010}

\begin{proof}
By \eqref{p1t}, the process $\{\hat p_t\}_{t=0}^T$ is $\{\mathcal G_t\}$-adapted. By Lemma \ref{p2t}, the conditional law of the next-step observable vector $Y_{t+1}=(x_{t+1}^\pi,l_{t+1})$ given the current information and the current control is completely determined by $(x_t^\pi,l_t,\hat p_t,u_t^\pi)$. Indeed, under regime $i\in\mathcal{M}$, the law of $Y_{t+1}$ is given by the kernel
$G_t^{(i)}(\cdot\,;x_t^\pi,l_t,u_t^\pi)$, while the next-step posterior belief is updated through
\[
\hat p_{t+1}
=
\Psi_t(\hat p_t,Y_{t+1};x_t^\pi,l_t,u_t^\pi),
\]
where $\Psi_t$ is defined by \eqref{psihatpxlu}. Therefore, conditional on $(x_t^\pi,l_t,\hat p_t,u_t^\pi)$, the conditional law of $(x_{t+1}^\pi,l_{t+1},\hat p_{t+1})$ is fully determined, and for any bounded measurable function $\varphi$,
\begin{align}\label{evarphi4var}
&\mathbb E\left[
\varphi(x_{t+1}^\pi,l_{t+1},\hat p_{t+1})
\;\middle|\;
x_t^\pi=x,l_t=l,\hat p_t=\hat p, u_t^\pi=u
\right]
=
\int_{\mathbb R^2}
\varphi \bigl(y_1,y_2,\Psi_t(\hat p,y;x,l,u)\bigr)\,
\widehat G_t(\mathrm{d}y;x,l,\hat p,u),
\end{align}
where $\hat G_t$ is given by \eqref{hatgxlp}. Thus the controlled transition kernel depends on the observed history through the belief state $(x_t^\pi,l_t,\hat p_t)$. For $t=T$, the terminal condition follows directly from the definition of the value function \eqref{termcon}, which does not depend on $\hat{p}$. For $t=T-1,T-2,\cdots,0$, conditioning on the current state $(x_t^\pi,l_t,\hat p_t)=(x,l,\hat p)$, by the Bellman's principle, we have
\begin{align*}
\hat J_t^*(x,l,\hat p;w)
& =
\min_{\pi_t}
\mathbb{E}\left[
\hat J_{t+1}^*(x_{t+1}^\pi,l_{t+1},\hat p_{t+1};w)
+\lambda\int_{\mathbb{R}}\pi_t(u)\ln\pi_t(u)\,\mathrm{d}u
\;\middle|\;
x_t^\pi=x, l_t=l, \hat p_t=\hat p
\right]\\
& = \min_{\pi_t} \int_{\mathbb R}
\left[
\mathbb E\left(
\hat J_{t+1}^*(x_{t+1}^\pi,l_{t+1},\hat p_{t+1};w)
\;\middle|\;
x_t^\pi=x, l_t=l, \hat p_t=\hat p, u_t^{\pi} = u \right)
+\lambda\ln\pi_t(u)
\right] \pi_t(u) \mathrm{d} u.
\end{align*}
Then, by \eqref{evarphi4var}, it further satisfies the following Bellman recursion
\begin{align}\label{hatjxl00}
\hat J_t^*(x,l,\hat p;w)
& =
\min_{\pi_t}
\int_{\mathbb R}
\Bigl(
\mathcal Q_t(u;x,l,\hat p; w)
+\lambda\ln\pi_t(u)
\Bigr) \pi_t(u) \, \mathrm{d} u,
\end{align}
where $\mathcal Q_t(u;x,l,\hat p; w)$ is defined in \eqref{qq22}.
To verify that the Gibbs minimizer is well defined, use the complete-information value as a lower bound on the partially observed value. Conditional on a belief, a complete-information investor can implement the same partially observed policy. Hence the partially observed continuation value is at least the belief-weighted complete-information value for the same $w$. Taking the conditional expectation over the next state gives a lower bound on $\mathcal Q_t$ with positive action curvature
\[
\sum_{i\in\mathcal M}p_i B_{1,t}^{(i)}\sum_{j\in\mathcal M}P_{ij}\bar F_{1,t+1}^{(j)},
\qquad (p_1,p_2)=(\hat p,1-\hat p).
\]
There are finitely many dates and regimes. Completing the square in this lower bound gives $\mathcal Q_t(u)\ge a u^2-C M$ for constants $a,C>0$ and $M=1+x^2+l^2+w^2$. A fixed standard Gaussian policy at subsequent dates gives the upper bound $\mathcal Q_t(u)\le C'(M+u^2)$ by the finite conditional second moments and the linear dynamics. These bounds apply at each step of the backward induction, starting with the terminal value.

Writing $q(u):=\mathcal Q_t(u;x,l,\hat p;w)$ and $Z:=\int_{\mathbb R}\exp[-q(v)/\lambda]\,\mathrm dv$, the two bounds imply $0<Z<\infty$ and $|\ln Z|\le C'' M$. For the Gibbs density, the upper bound on $q$ over $[-1,1]$ gives a lower bound on $Z$ of the form $c\exp(-C'' M)$. Together with the lower bound on $q$, this bounds the Gibbs tails by a Gaussian outside $|u|\le C_0\sqrt M$ for a sufficiently large constant $C_0$. Its second moment and the absolute entropy integral are therefore bounded by a constant times $M$. The linear dynamics then give finite state second moments and entropy expectations over the finite horizon. The measurable Bellman kernel and the normalization define an admissible Gibbs feedback.
Then, for any admissible control density $\pi_t$,
\begin{align*}
\int_{\mathbb R}\bigl(q(u)+\lambda\ln\pi_t(u)\bigr) \pi_t(u)\,\mathrm{d}u
&=
\int_{\mathbb R}
\left(
-\lambda\ln\pi_t^*(u)-\lambda\ln Z+\lambda\ln\pi_t(u)
\right)\pi_t(u) \mathrm{d}u \\
&=
-\lambda\ln Z
+\lambda\int_{\mathbb R}\pi_t(u)\ln\frac{\pi_t(u)}{\pi_t^*(u)}\,\mathrm{d}u
=
-\lambda\ln Z+\lambda\,\mathrm{KL}(\pi_t\|\pi_t^*),
\end{align*}
where $\mathrm{KL}(\pi_t\|\pi_t^*)$ denotes the Kullback-Leibler divergence. Since $\mathrm{KL}(\pi_t\|\pi_t^*)\ge 0$, with equality if and only if $\pi_t=\pi_t^*$ a.e., it follows that $\pi_t^*$ is the unique minimizer of \eqref{hatjxl00}. Therefore, there exists a unique optimal feedback control of Gibbs form, whose density function $\pi_t^*(u|x, l,\hat p; w)=\pi_t^*(u)$ satisfies \eqref{pittt}. Substituting the optimal control $\pi_t^*$ into \eqref{hatjxl00} yields \eqref{hatjt11}.

\end{proof}

\subsection{Proof of Theorem \ref{negtilpthm}}\label{appen1122}

\begin{proof}
Under the auxiliary law \eqref{tildeassetre}, the components of the joint gross-return vector may be dependent within a period. Their first, second and cross moments are the corresponding prior-weighted regime moments. These joint moments suffice for the quadratic recursion. We also substitute $\varepsilon$ with $\tilde{p}$, $\pi^c$ with $\pi^{sub}$, and $J$ with $\tilde{J}$. The no-learning benchmark value function is denoted by $\tilde J_t^{sub*}(x,l, \tilde p; w)$. For fixed $w\in\mathbb R$, its Bellman recursion is
\begin{align*}\label{eq:nl-bellman}
\tilde J_t^{\mathrm{sub}*}(x,l,\tilde p; w)
=
\min_{\pi_t}
\int_{\mathbb R}\pi_t(u)
\Big(
\mathbb E_{t, x, l, \tilde p}\big[
\tilde J_{t+1}^{\mathrm{sub}*}(x_{t+1}^{\pi},l_{t+1},\tilde p_{t+1}; w)
\,\big|\, u_t^{\pi}=u  \big]
+\lambda\ln\pi_t(u)
\Big)\mathrm du,
\end{align*}
where the conditional expectation uses the auxiliary transition law induced by \eqref{tildeassetre}. Completing the square as in the complete-information proof, with the auxiliary joint moments and the deterministic continuation prior, gives the feedback \eqref{pi316}. The prior mixture preserves positive definiteness of the asset second-moment matrix, so the action curvatures remain positive.

\end{proof}

\subsection{Proof of Theorem \ref{thm:policy_improvement_unified}}\label{thm:policy_improvement_unified00}

\begin{proof}
Let $q(u):= Q_t^{m, \pi}(u; \xi, w)$ and $Z:=\int_{\mathbb R} e^{-q(v)/\lambda} \mathrm{d} v$. Define $\rho^*(u):=\frac{e^{-q(u)/\lambda}}{Z}$, as in \eqref{definpitrho}.
Then, for any admissible control density $\rho(\cdot|\xi)$, we have
\begin{align*}
\int_{\mathbb R}
\bigl(q(u)+\lambda\ln\rho(u)\bigr) \rho(u)\,\mathrm du
&=
\int_{\mathbb R}
\left(
-\lambda\ln\rho^*(u)-\lambda\ln Z+\lambda\ln\rho(u)
\right) \rho(u) \mathrm du\\
&=
-\lambda\ln Z
+\lambda\int_{\mathbb R}\rho(u)\ln\frac{\rho(u)}{\rho^*(u)}\,\mathrm du
= -\lambda\ln Z +\lambda\,\mathrm{KL}(\rho\|\rho^*),
\end{align*}
where $\mathrm{KL}(\rho\|\rho^*)$ denotes the Kullback-Leibler divergence.
Since $\mathrm{KL}(\rho\|\rho^*)\ge 0$, with equality if and only if $\rho=\rho^*$ a.e., it follows that the unique minimizer of $\rho\mapsto
\int_{\mathbb R}
\bigl(q(u)+\lambda\ln\rho(u)\bigr) \rho(u)\,\mathrm du$
is precisely $\rho^*$.
Therefore,
\begin{equation}\label{eq:one_step_improvement}
(\mathcal T_t^{m, \pi_t'}J_{t+1}^{m, \pi})(\xi)
\le
(\mathcal T_t^{m, \pi_t}J_{t+1}^{m, \pi})(\xi)
=
J_t^{m, \pi}(\xi;w).
\end{equation}

We next show that the improved policy $\pi'$ improves the value function at every time. At the terminal time, both policies have the same payoff, so 
$J_T^{m, \pi'}(\xi;w)=J_T^{m, \pi}(\xi;w)=\Phi(\xi; w)$.
Now suppose that $J_{t+1}^{m, \pi'}(\xi;w)\le J_{t+1}^{m, \pi}(\xi;w)$ for all admissible states $\xi$. Since the evaluation operator is monotone in the continuation function, we have
\[
(\mathcal T_t^{m, \pi_t'}J_{t+1}^{m, \pi'})(\xi)
\le
(\mathcal T_t^{m, \pi_t'}J_{t+1}^{m, \pi})(\xi).
\]
Combining this with \eqref{eq:one_step_improvement} yields
\[
J_t^{m, \pi'}(\xi; w)
=
(\mathcal T_t^{m, \pi_t'}J_{t+1}^{m, \pi'})(\xi)
\le
(\mathcal T_t^{m, \pi_t'}J_{t+1}^{m, \pi})(\xi)
\le
(\mathcal T_t^{m, \pi_t}J_{t+1}^{m, \pi})(\xi)
=
J_t^{m, \pi}(\xi;w).
\]
By backward induction, the result follows for all $t=0,1,\cdots,T$.
If $\pi_t'(\cdot\mid \xi)\neq \pi_t(\cdot\mid \xi)$ on a set of positive Lebesgue measure, then $\mathrm{KL}(\pi_t(\cdot\mid \xi)\|\pi_t'(\cdot\mid \xi))>0$, and the one-step inequality in \eqref{eq:one_step_improvement} is strict at that state. Propagating this strict improvement backward through the Bellman recursion yields strict improvement at some current or earlier time-state pair.
\end{proof}

\subsection{Proof of Theorem \ref{thm:policy_convergence}}\label{thm:policy_convergence00}
\begin{proof}
Monotonicity follows by Theorem \ref{thm:policy_improvement_unified}. At $t=T-1$, every policy has the same terminal continuation value, so the first Gibbs update gives the minimizing density and value at $T-1$. Suppose that after $n$ updates the policy is optimal from $T-n$ onward. Its continuation value at $T-n$ is then $J_{T-n}^{m,*}$. The next Gibbs update uses this value and therefore gives the minimizing density at $T-n-1$. It also preserves optimality at all later dates, where the continuation values are already optimal. Backward induction gives the stated equality at each date. The Gibbs variational identity gives uniqueness of the minimizing density up to null sets at a fixed state.
\end{proof}

\end{document}